\documentclass[paper,twoside]{article}
\usepackage{amsfonts,amsmath}

\newtheorem{Definition}{Definition}
\newtheorem{Proposition}{Proposition}
\newtheorem{Theorem}{Theorem}
\newtheorem{Corollary}{Corollary}
\newtheorem{Lemma}{Lemma}

\newenvironment{Proof}{\hfill\par\noindent{\bf Proof:}}{\hfill\rule{2mm}{2mm}\par\noindent}
\newenvironment{Remark}{\vspace{0.1cm}\noindent{\bf Remark.}}{\vspace{0.1cm}}
\newenvironment{Example}{\vspace{0.1cm}\noindent{\bf Example.}}{\vspace{0.1cm}}

\def\vphi{\varphi}
\def\eps{\varepsilon}

\def\d{\mathrm{d}}
\def\C{\mathbb{C}}

\def\Scal{\mathrm{Scal}}
\def\Ric{\mathrm{Ric}}

\def\tr{\mathrm{tr}}

\def\R{\mathbb{R}}

\def\FF{\mathcal{F}}

\def\DD{\mathcal{D}}

\def\HH{\mathcal{H}}
\def\VV{\mathcal{V}}

\def\m{\mathrm{m}}
\def\n{\mathrm{n}}
\def\bsob{{\mathrm{H}}_0}
\def\sob{{\mathrm{H}}}

\def\diver{\mathrm{div}}
\def\grad{\mathrm{grad}}

\def\spec{\mathrm{spec}}

\begin{document}

\flushbottom

\begin{center} {\bf \Large $\mathbf{L^2}$-Homogenization of Heat Equations on Tubular Neighborhoods}
\end{center}

\vskip0.5cm

\begin{center} {
O. WITTICH \\\small\em Department of Mathematics and Computer Science, \\ Eindhoven University of Technology, P.O. Box 513, \\ 5600 MB Eindhoven, The Netherlands}
\end{center}

\vskip0.5cm

\begin{quote}
{\bf Abstract.}{\sl We consider the heat equation with Dirichlet boundary conditions on the tubular neighborhood of a closed Riemannian submanifold of a Riemannian manifold. We show that, as the tube radius decreases, the semigroup of a suitably rescaled and renormalized generator can be effectively described by a Hamiltonian on the submanifold with a potential that depends on the geometry of the submanifold and of the embedding.}
\end{quote}

\vskip0.5cm

\section{Introduction}

In this paper, we study the heat equation on tubular neighborhoods around closed submanifolds of Riemannian manifolds. Let $L$ be an $l$-dimensional closed Riemannian manifold isometrically embedded into an $m$-dimensional Riemannian manifold $M$, $l<m$. We denote by $L(\eps)$ the tubular neighborhood of $L$ of radius $\eps > 0$. We investigate the behaviour of the heat semigroup on $L(\eps)$ assuming Dirichlet boundary conditions on $\partial L(\eps)$. It turns out that for small values of $\eps$, the solution of the heat equation can be effectively described mainly by the solution of a heat equation with potential on the submanifold alone. The potential reflects geometric properties of the submanifold and of the embedding. A proper formulation of this statement requires some modifications of the initial problem. If the semigroups are suitably rescaled and renormalized, we obtain an asymptotic perturbation problem on a fixed tube. As $\eps$ tends to zero, we actually obtain convergence to a limit semigroup. Since the proper formulation of the problem already requires some analysis, the main results are not
stated before Section \ref{5}.

The motivation to study this asymptotic problem grew out of work about Brownian motion conditioned to tubes around submanifolds. We wanted to understand the results of \cite{SmoWeiWit:07} in a different way, exploiting also another method of proof which is more analytic than the proof for $M = \R^n$ given in \cite{SidSmoWeiWit:04} which uses stochastic differential equations. Conditioned Brownian motion is intimately connected with Brownian motion absorbed at the boundary of the tube. The generator of absorbed Brownian motion is the Dirichlet Laplacian. This, and how we use these results on the heat equation to prove an analogous result as in \cite{SidSmoWeiWit:04} about weak convergence of the path measures of conditioned Brownian motion will be explained in \cite{pap3}.

This paper is the first one out of two. Here, we identify the limit semigroup and show that the semigroups converge strongly in a suitable $L^2$-space. Regarding the application above, this is not sufficient since we will have to consider also the restriction of these heat equation solutions to the zero set $L\subset M$ and this is not well - defined without additional smoothness. Therefore, in the second paper \cite{pap2} we will prove some compactness property which will imply that the semigroups actually converge strongly in Sobolev spaces of arbitrarily large Sobolev index. Note that we have to exclude time $t=0$ in both of these statements because the limit $\eps$ to zero is not continuous with respect to the initial condition.

After some preliminaries about the Dirichlet Laplacian and the tube geometry in Sections \ref{2} and {3}, we construct in Section \ref{4} a perturbation problem for quadratic forms associated to the generators that we actually want to investigate. In Section \ref{5} we formulate the main results of this paper, Theorem  \ref{MainMosco} and \ref{koerzmain}.

The idea of the proof is motivated by the geometry of the tube. We compare the Laplacians of the induced metric on the tube with the Laplacians associated to a {\em reference metric} induced by the so called {\em Sasaki metric} on the normal bundle. If $L(1)$ is equipped with the reference metric, the projection $\pi :L(1)\to L$ turns out to be a Riemannian submersion with totally geodesic fibers. The Laplacian associated to this metric is therefore decomposable into a horizontal and a vertical part and allows for quite explicit calculations. On the other hand, the reference metric also provides the leading terms in the asymptotic expansion of the of the induced metric as $\eps$ tend to zero and both degenerate in the same way. That makes it natural to think of the one metric as a pertubation of the other. This point of view is developed by calculations in local coordinates in Section \ref{6}.

The limit dynamics takes place on a subspace of the full Hilbert space which can be identified with the kernel of the renormalized vertical operator introduced in Section \ref{7}. Here we also explain, why we may use the smallest eigenvalue of the Dirichlet problem for the flat unit ball $B\subset\R^{m-l}$ to renormalize the quadratic forms associated to the generators of the dynamics.

In Section \ref{8}, we finally collect some facts about epi-convergence and its relation to strong resolvent convergence and convergence of the associated semigroups. Then we prove epi-convergence for the quadratic forms associated to reference and induced metric and conclude from that Theorem \ref{MainMosco}.

\section{The Dirichlet laplacian on small tubes}\label{2}

We consider the {\em Dirichlet Laplacian} $\Delta_{\eps}$ on the Hilbert space $L^2(L(\eps),g)$ given by
$\Delta_{\eps} u := \Delta u$ for all $u$ in the domain $\DD (\Delta_{\eps}) :=\bsob^1\cap \sob^2 (L(\eps),g)$. $\Delta
:= - \diver_{g} \grad$ denotes the Laplace - Beltrami operator on $M$ and the definition is understood in terms of
Sobolev derivatives. By \cite{Tay:99}, Ch. 5.1, $\Delta_{\eps}$ is a positive and self adjoint operator with discrete, semi-simple spectrum for whom all eigenfunctions belong to $C^{\infty}(\overline{L(\eps)})$. Furthermore, by {\em Friedrich's construction} (\cite{Kat:80}, VI.3), it is the operator associated to the quadratic form
\begin{equation}
\label{unskaliert}
    q_{(\eps)} (u) := \int_{L(\eps)} \m (dp) \,\Vert \d u\Vert^2
\end{equation}
with domain $\DD (q_{\eps}) = \bsob^1 (L(\eps),g)$. Here, $\m$ denotes the {\em Riemannian volume measure} associated to $g$ and $\Vert - \Vert$ the norm induced by $g$.

\begin{Example} As a special case of this construction, we consider the case $L=\lbrace 0\rbrace \in\R^{m-l}$ and $\eps =1$. Then, $L(1)=B$ is the unit ball $B\subset\R^{m-l}$ and we denote the corresponding Dirichlet laplacian by $\Delta_B$, its spectral values by
\begin{equation}\label{balleval}
\spec \,\Delta_B =\lbrace 0<\lambda_0 < \lambda_1 < \lambda_2< ... \rbrace
\end{equation}
and the corresponding eigenprojections by $P_0,P_1,P_2,...$. Note that the eigen\-space associated to the smallest eigenvalue $\lambda_0$ is one-dimensional.
\end{Example}

\section{Tube geometry and reference metric}\label{3}

First we collect some basic facts about the tubular neighborhood of
$L\subset M$. Let $\varphi : L
\to M$ be an isometric embedding of the closed manifold $L$ into $M$. We consider the {\em normal bundle} $NL
\subset TM\vert_L$. The {\em exponential map} $\exp^M : TM \to M$ restricted to $NL$ thus yields a smooth map $\exp^{\perp} : NL
\to M$. By compactness of $L$, there is some $r > 0$, the {\em injectivity radius}, such that $\exp^{\perp} : U_{\eps}(0)\to L(\eps)$ is a diffeomorphism from
the open $\eps$-neighborhood $U_{\eps}(0)$ of the {\em zero section} in $NL$ to the {\em tubular $\eps$-neighborhood} of
$L\subset M$ for all $\eps < r$. In the sequel, we assume for simplicity that $r>1$. The {\em bundle projection} $\pi_N
: NL \to L$ induces a submersion $\pi =
\vphi\circ \pi_N \circ\exp^{\perp \,-1} = \exp^{\perp}\circ \pi_N \circ\exp^{\perp\,-1} : L(1) \to L$. For $q\in L$, the preimages $\pi^{-1}
(q)$ provide a decomposition of $L(1)$ into relatively closed, $m-l$-dimensional submanifolds. The decomposition $\FF :=
\cup_{q\in L} \pi^{-1} (q)$ forms a {\em codimension-$l$-foliation} of $L(1)$ in the sense of \cite{Ton:97}, Ch. 1. The
submanifolds $\pi^{-1}(q)$ are called the {\em leaves}.

On the tube, we consider the family of diffeomorphisms $\sigma_{\eps}:L(\eps)\to L(1)$ given by
\begin{equation}
\label{Reskalierung}
    \sigma_{\eps} (p) := \exp^{\perp}\circ ( \eps^{-1} \,\exp^{\perp\,-1} (p)),
\end{equation}
$\eps > 0$. These maps are called {\em rescaling maps}.

Apart from the metric $g$ induced from $M$ via the embedding, we will now construct another metric, the reference metric $g_0$ on the tubular neighborhood $L(\eps)$.
This is the metric induced by the {\em Sasaki metric} (cf. \cite{Sak:97}, (4.6), p. 55 and (4.11), p. 58). The covariant derivative on the normal bundle (cf. \cite{Cha:93}, 2.9. p. 89 f.) associated to the embedding
induces a connection map $K^{\perp}: TNL \to NL$, (cf. \cite{Sak:97}, p. 58) and a decomposition of the tangent spaces $T_{\xi}NL$, $\xi\in N_pL$ into {\em vertical} and {\em horizontal} subspaces
$T_{\xi}NL = \VV_{\xi} \oplus \HH_{\xi}$ where $\VV_{\xi} = \ker \pi_*=T_{\xi}N_pL$ and $\HH = \ker K^{\perp}$. Thus, every $\eta = \eta_V + \eta_H$ can be decomposed into its {\em vertical} and {\em horizontal} component. Furthermore, the tangent vector $\eta\in T_{\xi}NL$ can be identified with $(A,B)\in T_pL \oplus N_pL$ where
\begin{equation}\label{HVdeco}
\begin{array}{cc} A = \pi_*(\eta), & B = K^{\perp}(\eta). \end{array}
\end{equation}

\begin{Definition}\label{Referenzmetrik} Identify $\eta,\eta'\in T_{\xi}NL$ with $(A,B),(A',B')\in T_pL\oplus N_PL$ according to (\ref{HVdeco}). The {\em Sasaki metric} on $NL$ is given by
$$
\langle \eta,\eta'\rangle_0 = \langle A,A' \rangle_L + \langle B,B'\rangle_{NL}
$$
where the first scalar product is computed with respect to the Riemannian metric on $L$ and the second with respect to the bundle metric on $NL$. The Sasaki metric makes $NL$ a Riemannian manifold. Via $\exp^{\perp}$, the restriction of the Sasaki metric to $U_1(0)$ is transported to a Riemannian metric on $L(1)$ which we denote by $g_0$. In the sequel, $g_0$ will be referred to as the
{\em reference metric}.
\end{Definition}

\begin{Remark} {\rm (i)} If the tube is equipped with the reference metric, the projection $\pi : L(1)\to L$ turns out
to be a {\em Riemannian submersion} (meaning that $\pi_* :\HH_{\xi}\to T_{\pi(\xi)}L$ is an isometry for all $\xi\in TNL$) with {\em totally geodesic fibers} (cf. \cite{BerBou:82}, Example 2.5), a situation which is investigated in particular
in the third section of \cite{Vil:70}. Furthermore, note that with respect to the reference metric, the metric induced by the embedding of the leaves
$\pi^{-1}(q)$ makes the leaves isometric to the flat unit ball $B\subset \R^{m-l}$ and the metric induced on the zero section is the intrinsic metric $g_L$ on $L$.
{\rm (ii)} The decomposition $\eta = \eta_V + \eta_H$ of $\eta\in T_{\xi}NL$ is orthogonal with respect to $g_0$, i.e. $\langle\eta_V,\eta_H\rangle_0=0$.
{\rm (iii)} In contrast to the situation
for the induced metric, the Dirichlet Laplacian associated to the reference metric can be decomposed into two {\em
independent} parts, a {\em vertical} part along the leaves and a {\em horizontal} part along the submanifold. This is
proved in \cite{BerBou:82}, Theorem 1.5 for closed fibers. With slight modifications, the same holds for the boundary problem considered above. We will use this in the subsequent paper \cite{pap2}.
\end{Remark}

\section{Rescaling and renormalization}\label{4}

Denote the Riemannian volume measure associated to the reference metric by $\m_0$. First of all, $\m_0$ behaves under rescaling in the following way.

\begin{Lemma}\label{volumen} We have $\sigma_{\eps\,*}\m_0 = \eps^{l-m}\m_0$.
\end{Lemma}

\begin{Proof} The statement follows from the local coordinate expression (\ref{riemvoll}) below which is a consequence of Lemma \ref{Referenzmetrik lokal}.
\end{Proof}
Denote by
\begin{equation}
\label{Dichte}
    \rho := \frac{\d\m}{\d\m_0} \in C^{\infty}(\overline{L(1)})
\end{equation}
the {\em Radon-Nikodym density} of the volume measures associated to $g$, $g_0$, respectively. The density $\rho > 0$ is
strictly positive, hence pointwise multiplication of functions $f\in L^2 (L(\eps),g_0)$ by $\rho^{-1/2}$ yields a
unitary isomorphism between $L^2 (L(\eps),g_0)$ and $L^2 (L(\eps),g)$ for all $1\geq
\eps > 0$. This, together with the {\em rescaling map} from (\ref{Reskalierung}) yields a unitary isomorphism
$\Sigma_{\eps} : L^2(L(1),g_0)\rightarrow L^2(L(\eps),g)$ given by
\begin{equation}
\label{UniIso}
  \Sigma_{\eps}  u  := \frac{\sigma_{\eps}^*u}{\sqrt{\rho\eps^{m-l}}}
\end{equation}
where we use the shorthand $\sigma_{\eps}^*u := u\circ\sigma_{\eps}$. The maps $\Sigma_{\eps}$ restrict to isomorphisms
of the respective domains. To be precise, $\sqrt{\rho}\in C^{\infty}(\overline{L(1)})$ together with the fact that the
rescaling maps $\sigma_{\eps}$ are diffeomorphisms implies that for all $\eps > 0$, the maps
\begin{equation*}
\Sigma_{\eps}: \bsob^1 (L(1),g_0)\rightarrow
\bsob^1 (L(\eps),g)
\end{equation*}
are topological isomorphisms mapping the boundary Sobolev space $\bsob^1 (L(1),g_0)$ homeomorphically onto the domains of the quadratic forms $q_{\eps}$ associated to the respective Dirichlet operators.

To finally construct the perturbation problem that we will actually investigate, we will {\em rescale} now the quadratic forms (\ref{unskaliert}) to quadratic forms defined on $\bsob^1(L(1),g_0)$ and {\em renormalize} them with the help of the lowest eigenvalue $\lambda_0$ of the Dirichlet laplacian for the flat euclidean ball (\ref{balleval}).

\begin{Definition}\label{Die Reskalierten Funktionale} Let $\alpha > 0$, $\eps > 0$ and $\lambda_0 > 0$ be the smallest eigenvalue from (\ref{balleval}). The {\em
rescaled forms}
$$
F_{\eps,\alpha}(u) : \bsob^1 (L(1),g_0)\rightarrow \R
$$
are given by
\begin{equation}\label{resren}
F_{\eps,\alpha}(u) := \frac{1}{2}\left\lbrace q_{(\eps)}(\Sigma_{\eps}u) +(\alpha -
\frac{\lambda_0}{\eps^2})\langle\Sigma_{\eps}u,\Sigma_{\eps}u\rangle\right\rbrace
\end{equation}
where $\langle - , - \rangle$ denotes the scalar product on $L^2(L(1),g)$.
\end{Definition}

\begin{Remark} Note that the operators associated to the quadratic forms $F_{\eps,\alpha}$ are given by $\Delta (\eps) + \alpha$ where $\Delta(\eps)$ is given by
\begin{equation}\label{deltaeps}
\Delta (\eps) :=\Sigma_{\eps}^{-1}(\Delta_{\eps}  - \lambda_0/\eps^2)\Sigma_{\eps}
\end{equation}
with domain $\DD (\Delta(\eps)) := \bsob^1\cap\sob^2 (L(1),g_0)$. But $\Sigma_{\eps}: \bsob^1 (L(1),g_0)\rightarrow
\bsob^1 (L(\eps),g)$ is a unitary isomorphism and hence $\Delta (\eps)$ is also self-adjoint with discrete semi-simple spectrum
$$
\spec\,\Delta (\eps) = \left\lbrace \mu  -\frac{\lambda_0}{\eps^2}\,:\,\mu\in\spec\,\Delta_{\eps} \right\rbrace .
$$
\end{Remark}

The functionals $F_{\eps,\alpha}$, the associated operators $\Delta(\eps)$ and in particular their behavior as $\eps$ tends to zero will be what we will investigate in the remainder of this paper. The two main statements that we will prove in this paper will be formulated in the next subsection. As a preparation, we will first use partial integration to represent the functionals $F_{\eps,\alpha}$ in a slightly different way, thereby emphasizing the fact that they are actually functionals on a Hilbert space with weight given by the Riemannian volume of the reference metric.

\begin{Proposition}\label{modiDir1} The rescaled and renormalized form $F_{\eps,\alpha}$ is given by
\begin{equation*}
F_{\eps,\alpha}(u)= \frac{1}{2}\int_{L(1)} \m_0 (dp) \left\lbrace\Vert \d u
\Vert^2_{g(\eps)} + \left(W_{\eps}+\alpha - \frac{\lambda_0}{\eps^2}\right)u^2 \right\rbrace
\end{equation*} where $W_{\eps} := \sigma_{\eps}^{-1\,*}W$, $\Vert -\Vert_{g(\eps)}$ denotes the norm with respect to
the metric $g(\eps) :=\sigma_{\eps}^{-1\,*}g$ and the smooth potential $W$ is given by
\begin{equation*}
    W:= \frac{1}{2}\Delta \log\rho - \frac{1}{4}\Vert
    \d\log\rho\Vert^2\in C^{\infty} (\overline{L(1)}).
\end{equation*}
\end{Proposition}

\begin{Proof} Denote by $\langle -,-\rangle$ also the fiberwise scalar product associated to the induced metric $g$ and
by $\nabla$ the (vector field valued) gradient, i.e. $\nabla u:=\langle \d u,-\rangle$. We denote the scalar product for forms and vector fields by the same symbol. By
\begin{equation*}
    \d\Sigma_{\eps} u = \eps^{l-m/2}\,\rho^{-1/2} \,
\left(\d\sigma_{\eps}^*u
    -\frac{1}{2}\sigma_{\eps}^*u\,\,\, \d\log\rho \right),
\end{equation*}
we have
\begin{eqnarray*}
Q_{\eps}(u) &=& \eps^{l-m} \int_{L(\eps)} \m(dp)\,\rho^{-1} \left(\langle
\d\sigma_{\eps}^*u , \d\sigma_{\eps}^*u \rangle +
\frac{1}{4} \sigma_{\eps}^*(u)^2
\Vert\d\log\rho \Vert^2\right)  \\
& & - \frac{1}{2}\eps^{l-m} \int_{L(\eps)} \m (dp)\rho^{-1}
\langle\d\sigma_{\eps}^*(u)^2 , \d\log\rho\rangle \\
&=& \eps^{l-m} \int_{L(\eps)} \m(dp)\,\rho^{-1} \left(\langle
\nabla\sigma_{\eps}^*u , \nabla\sigma_{\eps}^*u \rangle +
\frac{1}{4} \sigma_{\eps}^*(u)^2
\Vert\nabla\log\rho \Vert^2\right)  \\
& & - \frac{1}{2}\eps^{l-m} \int_{L(\eps)} \m (dp)\rho^{-1}
\langle\nabla\sigma_{\eps}^*(u)^2 , \nabla\log\rho\rangle .
\end{eqnarray*}
Note now that by $\nabla\log\rho= -\rho\nabla\rho^{-1}$ the second term is given by
\begin{eqnarray*}
& & - \frac{1}{2}\eps^{l-m} \int_{L(\eps)} \m (dp)\rho^{-1}
\langle\nabla\sigma_{\eps}^*(u)^2 , \nabla\log\rho\rangle \\
&=& \frac{\eps^{l-m}}{2}\int_{L(\eps)}\m (dp)\langle\nabla\sigma_{\eps}^*(u)^2,\nabla\rho^{-1}\rangle =
\frac{\eps^{l-m}}{2}\int_{L(\eps)} \m (dp)\, \langle \nabla v, X \rangle
\end{eqnarray*}
where
\begin{equation*}
\begin{array}{cc}    v := \sigma_{\eps}^*(u)^2, & X:= \nabla \rho^ {-1}.
\end{array}
\end{equation*} By \cite{Tay:99}, Proposition 2.3, p. 128, we have for a smooth vector field $X$ and a smooth function
$u$ on the (possibly non - orientable) tube $L(\eps)$ that
\begin{equation*}
\int_{L(\eps)} \m(dp)\,  (v \,\diver X + \langle \nabla v, X\rangle) = \int_{\partial L(\eps)} \m_{\partial}(dp) \langle
X, \n\rangle
\end{equation*} where $\m_{
\partial}$ denotes the Riemannian volume on the boundary of the tube and $\n$ denotes the unit outward - pointing normal
to $
\partial L(\eps)$. Using the fact that we consider Dirichlet boundary conditions which implies vanishing of the boundary
terms, we obtain
\begin{equation*}
\frac{\eps^{l-m}}{2}\int_{L(\eps)}\m (dp)\langle\nabla\sigma_{\eps}^*(u)^2,\nabla\rho^{-1}\rangle  = -
\frac{\eps^{l-m}}{2}\int_{L(\eps)}\m (dp) \sigma_{\eps}^*(u)^2\,\diver(\nabla\rho^{-1}) .
\end{equation*} Finally, again by $\nabla\log\rho = -\rho \nabla\rho^{-1}$, and by $\diver\nabla u = \Delta u$, we
obtain
\begin{eqnarray*}
\diver(\nabla\rho^{-1}) &=& - \diver (\rho^{-1}\nabla\log\rho) =-\langle\nabla\rho^{-1},\nabla\log\rho \rangle -
\rho^{-1}\diver (\nabla\log\rho )\\
&=& \rho^{-1}\Vert\nabla\log\rho\Vert^2 - \rho^{-1}\Delta\log\rho.
\end{eqnarray*}
Adding up the different contributions yields
\begin{eqnarray*}
\int_{L(\eps)} \m (dp) \langle
\nabla \Sigma_{\eps}u ,\nabla \Sigma_{\eps}u\rangle &=& \int_{L(\eps)} \frac{\m (dp)}{\eps^{m-l}} \rho^{-1}(\langle
\nabla \sigma_{\eps}^*u ,\nabla \sigma_{\eps}^*u\rangle + W(\sigma_{\eps}^*u)^2 ) \\
&=& \int_{L(\eps)} \frac{\m (dp)}{\eps^{m-l}} \rho^{-1}(\Vert\d \sigma_{\eps}^*u\Vert^2 + W(\sigma_{\eps}^*u)^2 ).
\end{eqnarray*}
Then, by the transformation formula and Lemma \ref{volumen}
\begin{eqnarray*}
&& \int_{L(\eps)} \m (dp)
\rho^{-1}\eps^{l-m}\left(\Vert
\d \sigma_{\eps}^*u \Vert^2 + W \sigma_{\eps}^* (u)^2\right) \\
&=& \int_{L(1)} \eps^{l-m}\sigma_{\eps\,*}\m_0 (dp)\left(\Vert
\d u \Vert_{g(\eps)}^2 + \sigma_{\eps}^{-1\,*}W u^2\right)\\
&=& \int_{L(1)} \m_0 (dp) \left(\Vert \d u
\Vert_{g(\eps)}^2 + W_{\eps} u^2\right).
\end{eqnarray*}
That implies the statement.
\end{Proof}

\section{The homogenization results}\label{5}

It will turn out that the dynamics associated to the generators $\Delta
(\eps)$ approach a limit in $L^2 (L(1),g_{0})$
as $\eps$ tends to zero.

The cornerstone of the proof of the homogenization result is the following theorem about the
epi-limit of the rescaled functionals. It states that the limit functional remains bounded only on a subspace of the
whole domain. That results in the effect that the limit of the semigroups generated by $\Delta (\eps)$ is a semigroup on
the same subspace. This will be interpreted as homogenization along the fibers of the tube $L(1)$.

To state the first main result of this paper, we first have to introduce the subspace mentioned above. For that,
consider again the Dirichlet problem on the flat euclidean ball $B\subset \R^{m-l}$. The eigenspace to the lowest
eigenvalue $\lambda_0 > 0$ is one-dimensional and the eigenfunctions are orthogonally invariant. Let thus $U_0(\vert w
\vert)$ denote a normalized eigenfunction of the flat Dirichlet Laplacian with eigenvalue $\lambda_0$. On the tube, we
consider the function $u_0\in C^{\infty}(\overline{L(1)})$ given by
\begin{equation}
\label{vakuum} u_0(p) := U_0(d(p,\pi(p))) .
\end{equation}

\begin{Definition}\label{EE} Let $g_L$ denote the induced metric on the submanifold $L$. Then $E_{0}\subset L^2
(L(1),g_0)$ is given by
\begin{equation*} E_{0} := \lbrace u \in L^2(L(1),g_0) \,:\, u(p) = u_0(p)\, v(\pi(p)), v\in L^2(L,g_L)\rbrace .
\end{equation*}
\end{Definition}

{\bf Notation.} {\rm (i)} In the sequel, we will denote $E_{0}$ and the $L^2(L(1),g_{0})$-orthogonal projection onto it
by the same symbol. {\rm (ii)} For $f\in L^2(L,g_{L})$, we will frequently write $\overline{f}=f\circ\pi$ and use the
common notion {\em basic functions} for them.

Up to the renormalization, the operators $\Delta(\eps)$ are basically unitary transforms of the self-adjoint Dirichlet
Laplacians $\Delta_{\eps}$ on the tube $L(\eps)$. Hence, the operators $\Delta (\eps)$ are self-adjoint on the space
$L^2 (L(1),g_0)$. First of all, we need {\em equi-coercivity} of the sequence of functionals.

\begin{Proposition}\label{koerzmain} There are constants $\alpha_0,K > 0$ such that for all $\alpha \geq \alpha_0$ we have
$$
F_{\eps,\alpha}(u)\geq \frac{1}{4}q_0(u) $$
for all $\eps \leq 1/ 2K$. Here, $q_0$ denotes the quadratic form of the Dirichlet laplacian on $L(1)$ associated to the reference metric, i.e.
\begin{equation}\label{kuhnull}
q_0(u) = \int_{L(1)}\m_0(dq)\Vert \d u \Vert_{0}^2 .
\end{equation}
\end{Proposition}

\begin{Remark} Note that on the space $\bsob^1(L(1),g_0)$, $q_0(u)$ is equivalent to the 1-Sobolev norm by Poincar\'e inequality.
\end{Remark}

That implies convexity of the
$F_{\eps,\alpha}$ for the parameter values under consideration and together with the fact that the quadratic forms $F_{\eps,\alpha}$ define the self-adjoint operators $\Delta(\eps) + \alpha$, we may also conclude that $\mathrm{spec} (\Delta(\eps) + \alpha)\subseteq\lbrack 0,\infty)$. Thus, the
operators are uniformly bounded below and generate therefore strongly continuous semigroups on $L^2(L(1),g_{0})$.

Our major aim is to prove convergence of the semigroups associated to the operators $\Delta(\eps)$,  $\eps>0$ on the
Hilbert space $L^2(L(1),g_0)$.  The way we proceed is based on the following observations:

\begin{Theorem}\label{MainMosco} There is some $\alpha_{0} > 0$ such that for all $\alpha \geq \alpha_{0}$ the
functionals $F_{\eps,\alpha}:\bsob^1(L(1),g_0)\rightarrow\R$ epi-converge to $$
F_{\alpha}:\bsob^1 (L(1),g_0)\rightarrow \overline{\R}=\R\cup\lbrace\infty\rbrace $$
as $\eps\to 0$ with respect to the weak topology on $\bsob^1(L(1),g_{0})$. Here
\begin{equation*} F_{\alpha}(u) := \left\lbrace
\begin{array}{ll}
\frac{1}{2}\int_{L}
\m_L(dq)(\Vert\d v\Vert^2_{L} + (W_{L}+ \alpha) v^2) &,u = u_0\overline{v}\in E_{0} \\
\infty&,\mathrm{else}
\end{array}\right.
\end{equation*} where $\overline{v}=v\circ\pi$, $W_{L}:=W\vert_{L}$, $W$ is the potential from Proposition
\ref{modiDir1} and $u_0$ is the eigenfunction to $\lambda_0$ of the Dirichlet problem for the euclidean unit ball as in
(\ref{vakuum}).
\end{Theorem}

From epi-convergence of the functionals, we can conclude strong convergence of the resolvents and subsequently convergence of the semigroups.

\begin{Theorem}\label{MainSemi} {\rm (i)} Let $\alpha \geq \lambda_{0}+\inf_{q\in L}W(q)$. For every $u\in
L^2(L(1),g_0)$, we have $$
\lim_{\eps\to 0}(\Delta (\eps) + \alpha )^{-1} u = E_0 \,(\Delta_L + W_L + \alpha)^{-1}\, E_0 u $$
strongly in $L^2 (L(1),g_0)$. {\rm (ii)} For every $u\in L^2(L(1),g_0)$, we have $$
\lim_{\eps\to 0}e^{-\frac{t}{2}\Delta (\eps)}u  = E_0\, e^{-\frac{t}{2}(\Delta_L + W_L)} \, E_0 u $$
strongly in $L^2 (L(1),g_0)$ and uniformly for $t\in K$ where $K\subset (0,\infty)$ is compact.
{\rm (iii)} The statement of part {\rm (ii)} still holds true if we substitute $u$ by a sequence $u_{\eps}, \eps > 0$ with $\lim_{\eps\to 0} u_{\eps} = u\in L^2(L(1),g_0)$, i.e. we have $$
\lim_{\eps\to 0}e^{-\frac{t}{2}\Delta (\eps)}u_{\eps}  = E_0\, e^{-\frac{t}{2}(\Delta_L + W_L)} \, E_0 u $$
strongly in $L^2 (L(1),g_0)$ and uniformly for $t\in K$ where $K\subset (0,\infty)$ is compact.
\end{Theorem}

An expression for the potential $W_{L}$ which is the restriction of the potential $W$ from Proposition \ref{modiDir1} in terms of in- and extrinsic geometric quantities is given in Corollary \ref{Asymptotik Potential} in Section
\ref{EffPot}.

\section{An asymptotic formula for the metric}\label{6}

Recall that Proposition \ref{modiDir1} implies that the functionals $F_{\eps,\alpha}$ are understood up to the understanding of $g(\eps)$ and $W_{\eps}$ and that $W_{\eps}$ can be computed using the Radon-Nikodym density $\rho$. To understand $\rho$ and
$g(\eps)$ explicitly requires much more information about the geometry of the tube. This information will be collected now.

We have to analyze the geometry of the tubular
neighborhood more closely. We will investigate the asymptotic
behavior of the induced metric as the rescaling parameter $\eps$ tends to zero. It turns out that we can think of the
induced metric as a perturbation of the reference metric, and, even more important, the two metrics degenerate in the
same way.  Our investigation is thus based on the comparison of induced and reference metric as the tube radius tends to
zero.

The main work of this section is done in the proof of Theorem \ref{Fermi Metrik} where the induced metric coefficients
are calculated in suitable coordinates. In a series of corollaries, local expression for the dual metric, the Radon -
Nikodym density $\rho$ from (\ref{Dichte}) and the effective potential $W$ from Proposition \ref{modiDir1} are deduced
from this result without further complications.

The dynamical properties of the heat flow for the reference metric, namely those that follow from the fact that the
natural projection turns out to be a Riemannian submersion with totally geodesic fibers, will be discussed in the next
section.

\subsection{Fermi coordinates and the reference metric}

To analyze the geometry of the tubular neighborhood more
closely we will use suitably chosen local coordinates. They are obtained as follows:

Let $x = (x^i)_{i= 1,...,l}$ be arbitrary local coordinates for $L$ defined on an open subset $U\subset L$ with the
additional property that $NL\vert_U$ is trivial. This is the case for example if $U$ is contractible. Let
$\nu_{\alpha\,,\,\alpha = 1,...,m-l}$ denote {\em orthonormal sections} of $NL\vert_U$. The embedding of $L$ into $M$ is
still denoted by $\vphi : L\subset M$. Let now $w = (w^{\alpha})_{\alpha = 1,...,m-l}\in B:= \lbrace \vert w \vert < 1
\rbrace\subset\R^{m-l}$. Then the map $\tau (x,w) := \vphi(x) + w^{\alpha} \nu_{\alpha}(x)$ provides a local
trivialization of $NL\vert_U$ and
\begin{equation}
\label{Fermi}
    \Phi (x,w) := \eps^{\perp}\circ\tau (x,w) = \exp^M_{\vphi(x)} \left(w^{\alpha}
\nu_{\alpha}(x)\right)
\end{equation}
(summation convention) yields local coordinates for $L(1)\vert_U$, so called {\em Fermi coordinates}. Fermi coordinates
are {\em foliated} in the sense that the individual leaves are given by sets $\lbrace (x,w)\,:\,\vert w \vert < 1
\rbrace$ with fixed $x$. Note that we use small latin indices for the $x$-coordinates and small greek indices for the
$w$-coordinates in order to distinguish between parallel and transverse direction.

The rescaling maps $\sigma_{\eps}:L(\eps)\to L(1)$ defined in (\ref{Reskalierung}) are given in Fermi coordinates by $$
\sigma_{\eps} (x,w) = (x, w/\eps):U\times \lbrace \vert w \vert < \eps \rbrace \to U\times \lbrace \vert w \vert < 1
\rbrace. $$
Now we compute the reference metric in Fermi coordinates.

\begin{Lemma}\label{Referenzmetrik lokal} In local Fermi coordinates (cf. (\ref{Fermi})) we have
\begin{equation*}
g_{0}(x,w) = \left( \begin{array}{c|c} g_{L,ij} + w^{\mu}w^{\nu}C_{i\mu}^{\alpha}C_{j\nu}^{\beta}\delta{\alpha\beta} & w^{\mu}C_{i\mu}^{\alpha}\\\hline w^{\mu}C_{j\mu}^{\beta} & \delta_{\alpha\beta} \end{array}\right)
\end{equation*}
where $g_L$ denotes the metric on $L$and $C_{i\alpha}^{\mu} (x) :=
\langle \nu_{\mu}(x), D_i \nu_{\alpha}(x)\rangle$ denote the {\em connection coefficients} of the {\em induced
connection} $D$ on the normal bundle.
\end{Lemma}

\begin{Proof} By \cite{Sak:97}, p. 58, we have $\begin{array}{cc} K^{\perp} (\partial_i) = \partial_i + w^{\mu}C_{i\mu}^{\nu}\partial_{\nu}, & K^{\perp}(\partial_{\alpha}) = \partial_{\alpha}, \end{array}$
where $\partial_i = \Phi_* (\partial/\partial x^i)$, $\partial_{\alpha} = \Phi_* (\partial/\partial w^{\alpha})$. By Definition \ref{Referenzmetrik} that implies
$$
g_0(x,w) = \left( \begin{array}{c|c} 1 & w^{\mu}C_{i\mu}^{\beta}\delta_{\alpha\beta} \\\hline 0 & 1 \end{array}\right)  \left( \begin{array}{c|c} g_{ij}  & 0\\\hline 0 & \delta_{\alpha\beta} \end{array}\right) \left( \begin{array}{c|c} 1 & 0 \\\hline w^{\mu}C_{j\mu}^{\alpha}\delta_{\alpha\beta} & 1 \end{array}\right)
$$
which yields the statement above.
\end{Proof}

In particular, this result implies the remark {\rm (i)} above that with the metric induced by the reference metric, the leaves $(\pi^{-1}(q))$ are isometric to the flat unit ball
and that the zero section $(L,g_L)$ are isometrically embedded. This property corresponds to
the usual properties of normal coordinates and will allow us to compute the spectrum of the leaf operators which will be
introduced in the subsequent section. Furthermore, this result shows that in local Fermi coordinates, the Riemannian volume is given by
\begin{equation}\label{riemvoll}
\m_0(dxdw) = \sqrt{\det g_L} \,dx dw
\end{equation}
which is an alternative way to prove Lemma \ref{volumen}.

\subsection{The metric in local coordinates}

As explained above, the following proposition where the induced metric
on the tube is computed in local Fermi coordinates, is the basic result of this section.

\begin{Theorem}\label{Fermi Metrik} In local Fermi coordinates the metric tensor is given by
\begin{equation*}
    g(x,w) =
    \left(
\begin{array}{c|c}
    1 & cb^{-1} \\
\hline 0 & 1
    \end{array}\right) \,
    \left(
\begin{array}{c|c}
    a  & 0 \\
\hline 0 & b
    \end{array}\right)\,
    \left(
\begin{array}{c|c}
    1 & 0 \\
\hline b^{-1} c^+ & 1
    \end{array}\right)
\end{equation*} where
\begin{eqnarray*}
a_{ij} (x,w) &=& g_{L,ij} (x) - 2 w^{\alpha} g_{L,is} A_{\alpha\,s}^i
\\
& & +  w^{\alpha} w^{\beta} ( g_{L,rs}A_{\alpha\,i}^r A_{\beta\,j^s} - R_{i\alpha
j\beta}) + O(\vert w
\vert^3),
\\
c_{i\sigma} (x,w) &=& w^{\alpha} C_{i\alpha}^{\sigma} (x) + O(\vert w
\vert^2),
\\
b_{\mu\sigma} (x,w) &=& \delta_{\mu\sigma} -
\frac{1}{3}w^{\alpha}w^{\beta} R_{\mu\alpha\sigma\beta} (x) + O(\vert w \vert^3).
\\
\end{eqnarray*}
Here, $g_L$ denotes the metric on $L$, $A_{\alpha}$
the {\em Weingarten map} of the embedding, $C_{i\alpha}^{\mu} (x) :=
\langle \nu_{\mu}(x), D_i \nu_{\alpha}(x)\rangle$ denote the {\em connection coefficients} of the {\em induced
connection} $D$ on the normal bundle and $R$ is the {\em Riemannian curvature tensor} of $M$.
\end{Theorem}

\begin{Remark} For the local calculations, we use the Einstein summation convention in the sense that we sum over all indices occuring twice, even if
they are both lower or upper, i.e. the expression $R_{\mu\alpha\mu\beta}$ means $\sum_{\mu}R_{\mu\alpha\mu\beta}$ --
having still in mind that greek indices run from $1,...,m-l$ and small latin indices from $1,...,l$. Note in particular that by $g_{\alpha\beta}(x,0)=\delta_{\alpha\beta}$ the coefficients with {\em greek} indices are not affected by raising and lowering indices.
\end{Remark}

\begin{Proof} Let $\xi\in N_pL$ and $\eta\in T_{\xi}NL$ represented by $(A,B)\in T_pL\oplus N_pL$ as before. By \cite{Sak:97}, Lemma 4.8, p. 59., we have that
$$
Y(t) := D_{t\xi} \exp^{\perp} (A,tB)
$$
is a Jacobi vector field along the geodesic $\gamma$ given by $\gamma (0) = q$ and $\dot\gamma (0) = \xi$ with initial conditions $Y(0) = A$, $\nabla_t Y(0) = B + W_{\xi} A$ and $W_{\xi}$ denotes the {\em Weingarten map}. In particular, we obtain by
$$
Y (1) = D_{\xi} \exp^{\perp} (A,B) = D_{\xi} \exp^{\perp} (\eta)
$$
the desired tangent map. By Taylor expansion
$$
Y(t) = \tau_t \left(  Y(0) + t \, \nabla_tY(0) + \frac{t^2}{2} \,\nabla_t^2Y(0) + \frac{t^3}{6} \,\nabla_t^3Y(0) + O(t^4) \right)
$$
where $\tau_t : T_pM \to T_{\gamma(t)}M$ denotes parallel translation. Now clearly
$$
\begin{array}{lll} Y(0) = A, & \nabla_t Y(0) = B + W_{\xi} A, & \nabla_t^2 Y(0) = - R(\xi,A)\xi . \end{array}
$$
It remains to compute the third derivative. We have by the Jacobi field equations and the geodesic equation $\nabla_{\dot\gamma}\dot\gamma = 0$ in local coordinates
\begin{eqnarray*}
- \nabla_t^3 Y &=&  \nabla_t R(\dot\gamma,Y)\dot\gamma \\
&=& \nabla_t \left\lbrace R^K_{LMN} (\gamma(t))\dot{\gamma}^L\,Y^M \,\dot{\gamma}^N\partial_K \right\rbrace \\
&=&  \frac{d}{dt}\left\lbrace R^K_{LMN} (\gamma(t))\dot{\gamma}^L\,Y^M \,\dot{\gamma}^N \right\rbrace\,\partial_K +  R^K_{LMN} (\gamma(t))\dot{\gamma}^L\,Y^M \,\dot{\gamma}^N\,\nabla_t\partial_K .
\end{eqnarray*}
Hence at $t=0$, we obtain
$$
\nabla_t^3 Y (0) =  - R(\xi,B + W_{\xi}A)\xi + O(\vert\xi\vert^3) =  - R(\xi,B )\xi + O(\vert\xi\vert^3) .
$$
Inserting $t=1$ yields thus finally
$$
Y(1) = \tau_1 \left(  A + B + W_{\xi} A - \frac{1}{2} \,R(\xi,A)\xi - \frac{1}{6} \,R(\xi,B )\xi + O(\vert\xi\vert^3) \right).
$$
We first compute the scalar product of general vectors $\eta,\eta'\in T_{\xi}NL$. Since the values of these scalar products are not affected by parallel translation, we obtain using the fact that the {\em Weingarten map} maps the tangent space into itself
\begin{eqnarray*}
& & \langle D_{\xi} \exp^{\perp} (\eta),D_{\xi} \exp^{\perp} (\eta')\rangle \\ &=& \langle A,A' \rangle + \langle B, B' \rangle  + \langle A, W_{\xi} A' \rangle + \langle A', W_{\xi} A \rangle  + \langle W_{\xi}A, W_{\xi} A' \rangle \\  &&- (\langle R(\xi,B )\xi,B'\rangle  + \langle B, R(\xi,B' )\xi\rangle  + \langle R(\xi,B )\xi,A'\rangle  + \langle A, R(\xi,B' )\xi\rangle )/6 \\ & &  -  (\langle A, R(\xi,A')\xi\rangle  + \langle  R(\xi,A)\xi,A'\rangle )/2 + O(\vert\xi\vert^3).\\
\end{eqnarray*}
Note that $B$ might as well be $O(\vert\xi\vert)$ as in the case of $\partial_i$ which leads to absorption of the corresponding terms into the remainder.

Inserting now the respective decompositions of the coordinate vector fields
$$
\begin{array}{lll}
\eta = \partial_i: & A = \partial_i, & B = w^{\alpha}C_{i \alpha}^{\beta}\partial_{\beta}, \\
\eta = \partial_{\alpha}: & A = 0 & B = \partial_{\alpha}.
\end{array}
$$
yields for the induced metric on $L(1)$
\begin{eqnarray*}
g_{ij}(x,w)
&=& g_{L,ij} + w^{\mu}w^{\nu}C_{i\mu}^{\alpha}C_{j\nu}^{\beta}\delta_{\alpha\beta} - 2 w^{\mu}g_{L,is}A_{\mu\,i}^{s} \\
&&+ w^{\mu}w^{\nu}(g_{L,rs}A_{\mu\,i}^{r}A_{\nu\,j}^{s}   - \langle R(\partial_{\mu},\partial_j)\partial_{\nu},\partial_i\rangle ) + O(\vert w\vert^3),\\
g_{i\alpha}(x,w) &=& w^{\mu}C_{i\mu}^{\beta}\delta_{\alpha\beta} + O(\vert w\vert^3)\\
g_{\alpha\beta}(x,w)
&=& \delta_{\alpha\beta} -\frac{1}{3}\,w^{\mu}w^{\nu}\,\langle R(\partial_{\mu},\partial_{\alpha})\partial_{\nu},\partial_{\beta} \rangle   + O(\vert w\vert^3) .
\end{eqnarray*}
Note that we changed the sign when representing the Weingarten map by the matrix $A$, namely we write in local coordinates $W_{\alpha}X = - A_{\alpha j}^i X^j\partial_i$. For the components of the curvature tensor, we use the convention $\langle R(\partial_{\mu},\partial_{\alpha}),\partial_{\nu},\partial_{\beta}\rangle = R_{\mu\alpha\nu\beta}$.

Letting finally $b_{\mu\sigma} = g_{\mu\sigma}$, $c_{i\sigma} = g_{i\sigma}$ and $a_{ij} = g_{ij} - (cb^{-1}c^+)_{ij}$
yields the statement.
\end{Proof}

That implies for the dual metric

\begin{Corollary}\label{Metrik Kotangentialbundel} The dual metric is given by
\begin{equation*}
    g^{*}(x,w) = g_0^{*}(x,w) + \frac{1}{3}w^{\mu}w^{\nu}\left(
\begin{array}{c|c}
    O(\vert w\vert) & O(\vert w
\vert^2)\\
\hline
    O(\vert w \vert^2) &
    R_{\mu\alpha\nu\beta}+ O(\vert w
    \vert^3)
    \end{array}\right).
\end{equation*}
The expression $O(\vert w \vert^a)$ denotes smooth matrix valued functions such that all matrix
coefficients are of order $O(\vert w\vert^a)$ uniformly as $\vert w \vert\to 0$.
\end{Corollary}

\begin{Proof} By Lemma \ref{Referenzmetrik lokal}, the local expression for $g_0$ is given by
\begin{equation*}
g_0^* (x,w)  =g_0 (x,w)^{-1} =   \left( \begin{array}{c|c} 1 & 0 \\\hline -w^{\mu}C_{j\mu}^{\beta} & 1 \end{array}\right)\left( \begin{array}{c|c} g_{L,ij}  & 0\\\hline 0 & \delta_{\alpha\beta} \end{array}\right) \left( \begin{array}{c|c} 1 & -w^{\mu}C_{i\mu}^{\alpha} \\\hline 0 & 1 \end{array}\right)
\end{equation*}
Furthermore, from Proposition \ref{Fermi Metrik} we obtain
\begin{equation*}
    g^{*}(x,w) = g(x,w)^{-1} =
    \left(
\begin{array}{c|c}
    1 & 0 \\
\hline - b^{-1} c^+ & 1
    \end{array}\right) \,
    \left(
\begin{array}{c|c}
    a^{-1}  & 0 \\
\hline 0 & b^{-1}
    \end{array}\right)\,\left(
\begin{array}{c|c}
    1 & -cb^{-1} \\
\hline 0 & 1
    \end{array}\right).
\end{equation*} The statement follows by inserting the expressions for $a$, $b$ and $c$ from Proposition \ref{Fermi
Metrik}.
\end{Proof}

By this statement we may now finally compute the asymptotic difference between the induced and the reference metric.

\begin{Corollary}\label{Asymptotik Metrik} Let $g(\eps) := \sigma_{\eps}^{-1\,*}g$ be the {\em rescaled metric}. For the
dual metric we have $g^{*}(\eps) = g_{0}^*(\eps) + H^*(\eps)$ where
\begin{equation*}
    \sup_{p\in L(1)} \Vert H^*(\eps) - H^*(0)\Vert = O(\eps)
\end{equation*} as $\eps$ tends to zero. Furthermore, in local coordinates,
\begin{equation*}
    H^*(0)=\frac{1}{3}w^{\alpha} w^{\beta} \left(
\begin{array}{c|c}
    0& 0 \\
\hline 0 & R_{\mu\alpha\nu\beta}\end{array}\right).
\end{equation*}
\end{Corollary}

\begin{Proof} By Corollary \ref{Metrik Kotangentialbundel}, we have in local coordinates
\begin{eqnarray*}
& &g^{*}(\eps) - g_{0}^*(\eps) \\
&=& \left(\sigma_{\eps}^{-1\,*}(g - g_{0})(x,w)\right)^{-1}= \left(\begin{array}{c|c} 1 & 0 \\
\hline 0 & 1/\eps
\end{array}\right) (g^{*} - g_{0}^{*})(x,\eps\,w)\left(
\begin{array}{c|c} 1 & 0 \\
\hline 0 & 1/\eps
\end{array}\right) \\
&=&\frac{1}{3}w^{\alpha}
w^{\beta}\left(
\begin{array}{c|c} 0 & 0 \\
\hline 0 & R_{\mu\alpha\nu\beta}
\end{array}\right) + O(\eps).
\end{eqnarray*}
The statement follows now again from the fact that all matrix coefficients in the Taylor approximation are smooth
functions on the compact manifold $L$.
\end{Proof}

\subsection{Logarithmic derivatives and the effective potential}\label{EffPot}

In this subsections, we compute the asymptotic
expressions for the density $\rho$ from (\ref{Dichte}) and an expression
of the effective potential by ex- and intrinsic geometric quantities. In the sequel, we denote by $\tr$ the {\em trace}
of a matrix.

\begin{Corollary}\label{Log Dichte} In Fermi coordinates, the logarithmic Radon - Nikodym density $\log \rho$ is given
by
\begin{equation*}
\log\rho (x,w) = - w^{\alpha} \tr A_{\alpha} -\frac{1}{2}w^{\alpha}w^{\beta} \left( \tr A_{\alpha} A_{\beta} +
\frac{1}{3}R_{\mu\alpha\mu\beta} + R_{i\alpha i\beta} \right) + O(\vert w \vert^3) .
\end{equation*}
\end{Corollary}

\begin{Proof} By Proposition \ref{Fermi Metrik}
\begin{equation*}
    \det g = \det\left(
\begin{array}{c|c}
    a & 0 \\
\hline 0 & b
    \end{array}\right) = \det\left(
\begin{array}{c|c}
    g_L & 0 \\
\hline 0 & 1
    \end{array}\right)\left(
\begin{array}{c|c}
    g_L^{-1} a & 0 \\
\hline 0 & b
    \end{array}\right) = \det g_0 \det b\det g_L^{-1} a
\end{equation*} which implies by (\ref{riemvoll})
\begin{equation*}
    \log\rho (x,w) = \log\sqrt{\det g / \det g_0} = \frac{1}{2}\tr \left(\log (g_L^{-1} a)\, +
    \log (b)\right).
\end{equation*} Inserting the expressions for $a$, $b$ from Proposition \ref{Fermi Metrik} together with
\begin{eqnarray*}
    & &\log (1+w^{\alpha}M_{\alpha} + w^{\alpha}w^{\beta}M_{\alpha\beta}+ O(\vert w \vert^3))\\
    &=& w^{\alpha}M_{\alpha} +
    w^{\alpha}w^{\beta}(M_{\alpha\beta}-\frac{1}{2}M_{\alpha}M_{\beta}) + O(\vert w
    \vert^3).
\end{eqnarray*}
yields the statement.
\end{Proof}

Recall, that in local coordinates, the inverse rescaling map $\sigma_{\eps}^{-1} : L(\eps)\to L(1)$ is given by
$\sigma^{-1}_{\eps} (x,w) = (x,\eps\,w)$ and therefore $\sup_{p\in L(1)} d(\sigma_{\eps}^{-1} (p), \pi (p)) <
\eps$ which implies that $\sigma_{\eps}^{-1}$ converges uniformly to the projection map as $\eps$ tends to zero. For the
asymptotics of the logarithmic gradient of the density $\rho$ we obtain the following statement.

\begin{Lemma}\label{Asymptotik LogAbleitung} As $\eps$ tends to zero, we have
\begin{equation*}
    \sup_{p\in L(1)} \vert \rho\circ \sigma_{\eps}^{-1} -
    \rho\circ\pi\vert = O(\eps) .
\end{equation*} where $\rho\circ\pi = 1$, and
\begin{equation*}
    \sup_{p\in L(1)} \vert \d \log \rho\circ \sigma_{\eps}^{-1} -
    \d\log \rho\circ\pi\vert = O(\eps) .
\end{equation*} Furthermore $\d\log \rho\circ\pi = -\tau \circ\pi$, where $\tau$ denotes the {\em tension vector field}
of the embedding.
\end{Lemma}

\begin{Proof} Clearly $\rho\in C^{\infty} (\overline{L(1)})$ and $\rho > 0$. Furthermore, all coefficients in the
expansion of $\log \rho$ in Corollary \ref{Log Dichte} are smooth functions on $L$. Thus, by compactness of $L$, we
obtain in local coordinates
\begin{equation*}
    \log\rho\circ \sigma_{\eps}^{-1} = \log\rho(x,\eps\,w) = O(\eps),
\end{equation*} which implies $\rho\circ \sigma_{\eps}^{-1} = 1 + O(\eps)$, and
\begin{equation*}
    \d\log\rho\circ \sigma_{\eps}^{-1} = \d\log\rho(x,\eps\,w) = - \tr A_{\alpha} d w^{\alpha} +
    O(\eps),
\end{equation*} which implies the second statement since $\tau:=\tr A_{\alpha} dw^{\alpha}$ with {\em Weingarten map}
$A_{\alpha}$ is the local expression for the {\em tension field} of the embedding (cf.
\cite{Jos:98}, (8.1.14), p. 319).
\end{Proof}

As the final result of this subsection, we express the effective potential from Proposition \ref{modiDir1} in terms of
in- and extrinsic geometric quantities. Some of these quantities are non-standard, they will be defined first.

\begin{Definition}\label{parttraces} Let $p\in L$ and $e_{1,p},...,e_{l,p}$ be an orthonormal base of $T_{p}L\subset
T_{p}M$. Let $\Ric_M$, $\mathrm{R}_M$ be the Ricci- and the Riemannian curvature tensor of $M$. Then we denote by $$
\overline{\mathrm{Ric}}_{M/L} (p) :=  \sum_{i=1}^l \Ric_{M,p}(e_{i},e_{i}) $$
and $$
\overline{\mathrm{R}}_{M/L} (p):= \sum_{i,j=1}^l \langle e_{i},
\mathrm{R}_{M,p}(e_{i},e_{j})e_{j}\rangle $$
the partial traces of Ricci- and curvature tensor with respect to the tangent space of the submanifold.
\end{Definition}

The definition is independent of the special orthonormal base chosen. Now we derive the geometric expression for the
potential.

\begin{Corollary}\label{Asymptotik Potential} As $\eps$ tends to zero, we have for the potential
\begin{equation*}
    \sup_{p\in L(1)} \vert W\circ \sigma_{\eps}^{-1} -
    W_{L}\circ\pi\vert = O(\eps) .
\end{equation*} Furthermore, the effective potential on the submanifold is given by
\begin{equation*}
    W_{L} = \frac{1}{2} \Scal_L - \frac{1}{4} \vert\tau\vert^2
    - \frac{1}{6} (\Scal_M + \overline{\Ric}_{M/L}
    + \overline{\mathrm{R}}_{M/L})
\end{equation*} where $\tau$ is the tension vector field of the embedding, $\overline{Ric}_{M/L}$ and
$\overline{\mathrm{R}}_{M/L}$ denote the partial traces of Ricci- and curvature tensor defined above, and $\Scal_M$,
$\Scal_L$ denote the respective scalar curvatures.
\end{Corollary}

\begin{Proof} Recall the description of the effective potential $W$ in Proposition
\ref{modiDir1}. By Corollary \ref{Log Dichte} and Lemma
\ref{Asymptotik LogAbleitung}, we have in local coordinates
\begin{eqnarray*}
\Delta \log\rho &=& \frac{1}{\sqrt{\det g}}\,
\partial_K (g^{KL}\sqrt{\det g}\,
\partial_L\log\rho)\\
&=& g^{KL}
\partial_{KL}\log\rho + (
\partial_K g^{KL})
\partial_L\log\rho + g^{KL}\frac{
\partial_K\sqrt{\det g}}{\sqrt{\det g}}\,
\partial_L\log\rho\\
&=& g^{\alpha\beta}(
\partial_{\alpha\beta}\log\rho +
\frac{
\partial_{\alpha}\sqrt{\det g}}{\sqrt{\det g}}\,
\partial_{\beta}\log\rho)(0) + O(\vert w\vert)\\
&=& g^{\alpha\beta}(
\partial_{\alpha\beta}\log\rho +
\frac{
\partial_{\alpha}\rho\sqrt{\det g_L}}{\rho\sqrt{\det g_L}}\,
\partial_{\beta}\log\rho)(0) + O(\vert w\vert)\\
&=& g^{\alpha\beta}(
\partial_{\alpha\beta}\log\rho +
\partial_{\alpha}\log\rho\,
\partial_{\beta}\log\rho)(0) + O(\vert w\vert) \\
&=& \tr (A_{\alpha})^2 - \tr (A_{\alpha}^2) -
\frac{1}{3}R_{\mu\alpha\mu\alpha} - R_{i\alpha i\alpha} + O(\vert w \vert)
\end{eqnarray*}
where the convention that we sum over all indices occuring twice is still in force. By {\em Gauss equation} (see
\cite{KobNom:69}, Prop. 4.1, p. 23
\begin{equation*}
    \tr (A_{\alpha})^2 - \tr (A_{\alpha}^2) = \Scal_L - R_{ijij}
\end{equation*} we obtain
\begin{eqnarray*}
\frac{1}{2}\Delta \log\rho &=&
\frac{1}{2}\left(\tr(A_{\alpha})^2-\tr (A_{\alpha}^2) -
\frac{1}{3}R_{\mu\alpha\mu\alpha} - R_{i\alpha i\alpha}\right) + O(\vert w
\vert)\\
&=& \frac{1}{2}\left(\Scal_L - \frac{1}{3}(R_{\mu\alpha\mu\alpha} + 3R_{i\alpha i\alpha} + 3R_{ijij})\right) + O(\vert w
\vert) \\
&=& \frac{1}{2}\Scal_L - \frac{1}{6}(\Scal_M +
\overline{\Ric}_{M/L} + \overline{\mathrm{R}}_{M/L}) + O(\vert w \vert)\\
\end{eqnarray*}
since $R_{\mu\alpha\mu\alpha} + 2R_{i\alpha i\alpha} + R_{ijij}=\Scal_M$, $R_{i\alpha i\alpha} + R_{ijij}=\Ric_{ii}=
\tr_{TL} \Ric_M$, and $R_{ijij}=\tr_{TL}\mathrm{R}_M$. Finally, $$
\Vert \d\log\rho\Vert^2\circ\pi = \vert\tau\vert^2\circ\pi$$
follows directly from Lemma \ref{Asymptotik LogAbleitung}.
\end{Proof}

\section{The vertical operator}\label{7}

Due to the renormalization, for the proof of epi-convergence we have to use also some spectral properties of what we will for now call the {\em vertical operator}. Let $\pi:L(1)\rightarrow L$ denote again the projection map. By Lemma \ref{Referenzmetrik lokal}, the fibers $\pi^{-1}(q)\subset L(1)$, equipped with the metric induced by $g_0$, are all Riemannian manifolds isometric to the euclidean unit ball $B\subset\R^{m-l}$. Thus, we may consider the associated Laplace-Beltrami operators $D_q$ on the fiber $\pi^{-1}(q)$.

\begin{Definition}\label{horivert} Let $f\in C^2 (L(1))$. The operator $$
\Delta_{0}^Vf (p):= D_{\pi (p)}\left(f\vert_{\pi^{-1}(\pi(p))}\right)(p)$$
is called the {\em vertical operator}.
\end{Definition}

\begin{Remark} In the local Fermi coordinates considered above, the vertical operator is given
by $\Delta_{0}^V = -\delta^{\alpha\beta}\,\partial^2/\partial w^{\alpha}\partial w^{\beta}$.
\end{Remark}

In order to make use of the preceding discussion for the investigation of the spectral properties of the vertical operator as an operator on $L^2(L(1),g_{0})$, we consider the tube $L(1)$ equipped with the reference metric $g_0$.
Note first that due to (\ref{riemvoll}), the Hilbert space with respect to $g_0$ obeys the {\em direct integral decomposition} (\cite{ReeSim:78})
$$
L^2(L(1),g_0) = \int^{\oplus}_L \m_L(dq) \,L^2(\pi^{-1}(q)) .
$$
Now we will construct an associated decomposable
self-adjoint operator which extends the vertical operator defined for functions $f\in C(\overline{L(1)})\cap C^2(L(1))$
with Dirichlet boundary conditions $f\vert_{\partial L(1)} = 0$.

Since the differential expression for the Laplacian depends smoothly on the fibers, the family is
measurable. All fibers can be isometrically mapped to the Euclidean unit ball $B\subset
\R^{m-l}$. Hence, the operators $D_q$ are unitarily equivalent to the Dirichlet Laplacian on $B\subset \R^{m-l}$ for all
$q\in L$. Hence, they are self-adjoint, their spectra $\spec (D_q)$ are {\em semi-simple} (i.e. consist only of
eigenvalues of finite multiplicity), and they are strictly positive with smallest eigenvalue $\lambda_0 > 0$. Hence, by
\cite{ReeSim:78}, Theorem XIII.85, p. 284, the operator
\begin{equation}
\Delta_{0}^V := \int_L^{\oplus} \m_L (dq) \, D_q
\end{equation}
is self-adjoint and extends the operator $\Delta_{0}^V$ defined on functions $f\in C(\overline{L(1)})\cap C^2(L(1))$
with Dirichlet boundary conditions $f\vert_{
\partial L(1)} = 0$. It is therefore justified to denote it by the same symbol. Again by \cite{ReeSim:78}, Theorem XIII.85, p. 284, the renormalized vertical operator
\begin{equation}
\label{DLO} \Delta_{0}^V - \lambda_0 := \int_L^{\oplus} \m_L (dq) \, (D_q -
\lambda_0)
\end{equation}
is as well self-adjoint and $$
\spec (\Delta_{0}^V-\lambda_{0}) := \lbrace \lambda_{k}-\lambda_{0}\,:\, k\geq 1\rbrace$$
where $\lambda_{k,k\geq 0}$ denotes the spectrum of the Dirichlet Laplacian on the Euclidean unit ball
$B\subset\R^{m-l}$. The eigenspace of the Dirichlet Laplacian in $B\subset\R^{m-l}$
belonging to the smallest eigenvalue $\lambda_0$ is one-dimensional, and the eigenfunctions are orthogonally invariant.
Let again $u_0 (w) := U_0 (\vert w
\vert)$ be a normalized eigenfunction generating this eigenspace. Recall, \cite{ReeSim:78}, XIII.16, Definition, p. 281,
that a bounded operator is called {\em decomposable}, if it can be written as a direct integral of operators as above.
Then the kernel coincides with the asymptotic subspace from Definition
\ref{EE}.

\begin{Corollary}\label{kern Dirichlet} {\rm (i)} The Dirichlet leaf operator $\Delta_{0}^V - \lambda_0$ is non-negative
with kernel
\begin{equation*}
\ker (\Delta_{0}^V - \lambda_0) = \lbrace u\in V\, :\, u(p):= u_0 (p)\, \overline{v}(p), v\in L^2 (L,g_L)\rbrace.
\end{equation*} Here, $u_0 (p) := U_0(d(p,\pi(p)))$, where $d$ denotes the Riemannian distance on $L(1)$ and $\overline{v}(p)=v(\pi (p))$. {\rm (ii)} The orthogonal projection $E_{0}$ onto $\ker (\Delta_{0}^V - \lambda_0)$ is decomposable.
\end{Corollary}

\begin{Remark} {\rm (i)} Note that in local Fermi coordinates $p\equiv (x,w)$ we have $$
u_0(p) = U_{0}(d(p,\pi(p))) = U_0(\vert w\vert) = u_0(w) $$
which justifies our slight abuse of notation. {\rm (ii)} In particular, the corollary implies that the kernel coincides
with the asymptotic subspace $E_0$ from Definition \ref{EE}.
\end{Remark}

\begin{Proof} {\rm(i)} A family $u_{q,q\in L}$ is contained in the kernel $\ker (\Delta_{0}^V -
\lambda_0)$ if and only if $(D_q - \lambda_0) u_q = 0$ for $\m_L$-almost all $q\in L$. The projection $E_{0,q}$ onto the
kernel $\ker (D_q - \lambda_0)$ is given by the $u_0$-weighted mean on the fiber, i.e.
\begin{equation*} E_{0,q} u_q (p) = u_0 (p) \,\langle u_0, u_q\rangle_q .
\end{equation*} Hence, $u(p) := u_{\pi(p)}(p)$ is contained in the kernel of $\Delta_{0}^V -\lambda_0$ if and only if
\begin{equation*} u(p) = E_{0,\pi(p)} u_{\pi(p)} (p) = u_0(p)\,\langle u_0, u_{\pi (p)}\rangle_{\pi (p)}
\end{equation*} where the equation is understood in terms of equivalence classes in $L^2(L(1),g_0)$. By the
Cauchy-Schwarz inequality, the function
\begin{equation*} v(q) := \langle u_0, u_{q}\rangle_{q}
\end{equation*} is contained in $L^2 (L,g_L)$ and $u (p) = u_0(p) \, v(\pi (p))$. {\rm(ii)} The projection onto the
kernel is given by
\begin{equation*} E_{0} = \int_L^{\oplus} \m_L (dq) E_{0,q} .
\end{equation*}
\end{Proof}

A fortiori, by \cite{ReeSim:78}, Theorem XIII.85, p. 284, we may conclude analogous results for all eigenspaces. We thus
obtain the following spectral decomposition of the operator $\Delta_{0}^V$: Let $\lambda_{k,k\geq 0} $ be the collection
of eigenvalues of the Dirichlet Laplacian on the euclidean unit ball $B\subset \R^{m-l}$ as above and $P_{k,k\geq 1}$
together with $P_{0}=E_{0}$ denote the corresponding eigenspaces. Again, we denote the eigenspaces and the corresponding
orthogonal projections by the same symbol. Then, by mapping the euclidean unit ball $B$ isometrically onto the fiber
$\pi^{-1}(q)$, we obtain projections $E_{k,q}$ on $L^2(\pi^{-1}(q),\m_{q})$ induced by $P_{k}$. Then,  by
\cite{ReeSim:78}, Theorem XIII.85, p. 284, we can compute the spectral decomposition of $\Delta_0^V$.

\begin{Proposition}\label{Dspektral} The operator $\Delta_{0}^V$ is self - adjoint on $L^2(L(1),g_0)$ with spectral
decomposition
\begin{equation*}\label{Dspektrum} \Delta_{0}^V = \lambda_{0} E_{0} + \sum_{k\geq 1}\lambda_{k} E_{k}
\end{equation*} where
\begin{equation*} E_{k} = \int_L^{\oplus} \m_L (dq) E_{k,q}.
\end{equation*}
\end{Proposition}

Finally, note that there is an analogous notion of direct integral decomposition for quadratic forms, and that the
quadratic form associated to the vertical operator is the direct integral of the respective quadratic forms on the
leaves given by
\begin{equation}\label{formvertikal}
q_V (u) := \int_{L}^{\oplus} \m_L (dq) \, (\Vert \d_q u_q
\Vert^2_q - \lambda_0 \Vert u_q \Vert_q^2)
\end{equation}
with domain
$$
\DD (q_V) := \int_{L}^{\oplus} \m_L (dq) \, \bsob^1 (\pi^{-1}(q)) .
$$
Here $\d_q$ denotes the exterior derivative on the leaf $\pi^{-1} (q)$, the leaf is equipped with the metric induced by the embedding, and $\Vert-\Vert_q$ denotes the norms on $\bsob^1 (\pi^{-1}(q))$ and $L^2 (\pi^{-1}(q))$, respectively.

\begin{Remark} Note that the spectral eigenspaces $E_{k}$ are infinite-dimensional in general. In the next subsection,
we will show that they all consist of eigenfunctions of the Dirichlet Laplacian $\Delta_{0}$ associated to $g_0$.
\end{Remark}

\section{Epi-convergence and convergence of the semigroups}\label{8}

In this section, we begin to discuss the case of the induced metric. We prove first
Theorem
\ref{epiMAIN} and conclude
from it the first part of our main result Theorem \ref{MainSemi} on the strong convergence of the family of semigroups
in the space $L^2(L(1),g_0)$. This is achieved by first proving strong resolvent convergence of the corresponding family of generators. Since the resolvents solve the Euler equations of a minimization problem associated to $F_{\eps,\alpha}$, they represent the unique minimizers of these variational problems. This is explained in detail in the proof of Theorem \ref{MainSemi}. Hence, convergence of these resolvents in $L^2(L(1),g_0)$ can be concluded from convergence
of the minimizers of the functionals $F_{\eps,\alpha}$ with respect to the weak topology in $\bsob^1(L(1),g_0)$ to the minimizer of the corresponding functional $F_{\alpha}$ as $\eps$ tends to zero. The main tool to prove this is {\em epi-convergence}. This is due
to the following fact from \cite{DaM:92}, Corollary 7.24, p. 84 which we review in a simplified form:

\begin{Proposition}\label{thereason} Suppose that the sequence of functionals $F_{n,n\geq
1}:X\to\overline{\R}=\R\cup\lbrace\infty\rbrace$ is equi-coercive and epi-converges to a function $F:X\to\overline{\R}$
with a unique minimizer $x\in X$. Let, for every $n\geq 1$, $x_n$ denote a minimizer of $F_n$. Then $x_n$ converges to
$x$ and $F_n(x_n)$ converges to $F(x)$ as $n$ tends to infinity.
\end{Proposition}

Hence, we will first prove equi-coercivity and epi-convergence of the functionals $F_{\eps,\alpha}$ to $F_{\alpha}$ with
respect to the weak topology of $\bsob^1(L(1),g_0)$. That implies convergence of the minimizers in this topology. By
another topological argument using the {\em Sobolev embedding theorem}, we can conclude from this strong convergence of
the resolvents in $L^2(L(1),g_0)$.

Before we come to the proofs, we summarize some definitions and facts about epi-convergence for the convenience of the
reader. They can all be found in the textbook \cite{DaM:92}. We begin with the notion of equi-coercivity.

\begin{Definition}\label{equikoerziv} A sequence $F_n : X \to \R\cup \lbrace\infty\rbrace$ of functions on a topological
space $X$ is called {\em equi-coercive}, if for every $t\in\R$ there exists a closed and countably compact subset
$K_t\subset X$ such that $\lbrace F\leq t\rbrace \subset K_t$.
\end{Definition}

According to Proposition 8.16 in \cite{DaM:92}, p. 97, and  Proposition 8.1, p. 86, epi-convergence in the weak topology of reflexive Banach spaces can be characterized as follows.

\begin{Definition} \label{epidef} Let $F_{n,n\geq 0}$ be a sequence of functionals on the reflexive Banach space $X$,
equi-coercive with respect to the weak topology. Then, $F_n$ {\em epi-converges} to the functional $F$ on $X$ {with
respect to the weak topology}, iff
\begin{enumerate}
\item For all $u\in X$ there is a weakly convergent sequence $u_n\to u$ such that
\begin{equation}
\lim_{n}  F_n(u_{n}) = F (u).
\end{equation}
\item For all $u\in X$ and for all weakly convergent sequences $u_n\to u$ we have
\begin{equation}
\liminf_{n}  F_n (u_{n}) \geq F (u).
\end{equation}
\end{enumerate}
\end{Definition}

\begin{Remark} Let the equi-coercive sequence $F_n(u)\geq \Vert u \Vert_X$ epi-converge to $F$ with respect to the weak
topology of the reflexive Banach space $X$, and let $\phi \in X^*$ be a linear functional. Then, by $$
F_n(u)+\phi(u) \geq  F_n (u) - \Vert \phi\Vert_{X^*}\,\Vert u \Vert_X \geq \frac{1}{2}\left(\Vert u \Vert_X^2 - \Vert
\phi\Vert_{X^*}^2 \right) $$
the sequence $F_n + \phi$ is also equi-coercive and since weak convergence $u_n\to u$ implies by definition that $\phi
(u_n)$ converges to $\phi(u)$, $F_n + \phi$ epi-converges to $F + \phi$ with respect to the weak topology on $X$. We
will use this fact freely in the sequel.
\end{Remark}

\subsection{Equi-coercivity of the reference family}

The basic idea to investigate the problem is to think of
the family $g(\eps)$ of rescaled induced metrics as a
perturbation of the family $g_0(\eps)$. Both degenerate in the same way, their difference tends to $H^*(0)$ as $\eps$ tends
to zero. Since the Laplacian is a nonlinear functional of the metric, we have to estimate the effect of the perturbation
$g(\eps) - g_0(\eps)$ on the final result. We will do this by comparing the functionals associated to the different
families. Therefore, we first prove some properties of the unperturbed family of functionals $F^0_{\eps,\alpha}$.

First of all, we will show that the rescaled an renormalized family $F^0_{\eps,\alpha}$ is uniformly
bounded below by the norm on $\bsob^1 (L(1),g_{0}))$ provided $\alpha > 0$ is large enough.

As a preparation, we have to clarify the relation of 1-Sobolev norm, horizontal and vertical operator. The quadratic form $q_0$ defining the Dirichlet laplacian on the tube associated to the metric $g_0$ can be decomposed according to
\begin{equation}\label{zerlege}
q_0(u) = \int_{L(1)}\m_0(dq) \Vert \d u\Vert_0^2 = q_V (u) + \tau (u) + \lambda_0 \langle u,u\rangle_0
\end{equation}
where $q_V$ is the form associated to the vertical operator (see (\ref{formvertikal})). We establish first some properties of the quadratic form $\tau$.

\begin{Lemma}\label{lokalfnull} In local Fermi-coordinates, we have
$$
\tau (u) = \int \m_0(dxdw)\,g_L^{ij} (\partial_iu - \frac{1}{2}C_{i\mu}^{\alpha}L^{\mu}_{\alpha}u)(\partial_ju - \frac{1}{2}C_{j\nu}^{\beta}L^{\nu}_{\beta}u)
$$
where $L^{\alpha}_{\mu} = w^{\alpha}\partial_{\mu} - w^{\mu}\partial_{\alpha}$ are Killing vector fields on the fibers. In particular
\begin{enumerate}
\item $\tau \geq 0$ is non-negative,
\item under rescaling, $\tau_{\eps}=\tau\circ\Sigma_{\eps}$ does not depend on the parameter $\eps > 0$.
\end{enumerate}
\end{Lemma}

\begin{Proof} By Lemma \ref{Referenzmetrik lokal}, we have in local coordinates
$$
\sigma_{\eps  \,*} \langle \d u , \d u \rangle_0 = \frac{1}{\eps^2}\delta^{\alpha \beta}\partial_{\alpha} u \, \partial_{\beta}u + g_L^{ij} (\partial_iu - w^{\mu}C_{i\mu}^{\alpha}\partial_{\alpha}u)(\partial_ju - w^{\nu}C_{j\nu}^{\beta}\partial_{\beta}u),
$$
hence $\tau$ is given by
$$
\tau (u) = \int_{U\times B(1)} dxdw\,\sqrt{\det g_L} \, g_L^{ij} (\partial_iu - w^{\mu}C_{i\mu}^{\alpha}\partial_{\alpha}u)(\partial_ju - w^{\nu}C_{j\nu}^{\beta}\partial_{\beta}u).
$$
Furthermore since $D$ is a {\em metric connection} which implies $C_{i\alpha}^{\mu} = - C_{i\mu}^{\alpha}$, we obtain
$$
w^{\mu}C_{i\mu}^{\alpha}\partial_{\alpha} = \frac{1}{2}(w^{\mu}C_{i\mu}^{\alpha}\partial_{\alpha} + w^{\alpha}C_{i\alpha}^{\mu}\partial_{\mu}) =  \frac{1}{2}C_{i\mu}^{\alpha}(w^{\mu}\partial_{\alpha} - w^{\alpha}\partial_{\mu}).$$
These vector fields generate orthogonal transformations of the fiber and therefore isometries. {\rm (i)} follows from the positive definiteness of $g_L$. For {\rm(ii)}, note first that in the rescaling map we have $\rho \equiv 1$ and thus $\Sigma_{\eps}u = \sqrt{\eps^{l-m}} \,\sigma_{\eps}^{*}u$. Hence in local coordinates
$$
\left( \partial_i -\frac{1}{2}C_{i\mu}^{\nu}L_{\mu}^{\nu}\right) (\Sigma_{\eps}u) = \eps^{\frac{l-m}{2}}\,\left( \partial_iu -\frac{1}{2\eps}C_{i\mu}^{\nu}L_{\mu}^{\nu}u\right) (x,w/\eps)
$$
and integration using Lemma \ref{volumen} implies the statement.
\end{Proof}

Equi-coercivity of the reference family is now given by the following statement.

\begin{Lemma}\label{koerziv} Let $\alpha  \geq \lambda_{0}$. Then we have for all $u\in\bsob^1 (L(1),g_{0})$ $$
F^0_{\eps,\alpha}(u) \geq \frac{1}{2}q_0(u). $$
\end{Lemma}

\begin{Proof} By Lemma \ref{lokalfnull}, we have
\begin{equation*}
F^0_{\eps,\alpha}(u) = \frac{1}{2}\left\lbrace\eps^{-2} q_V(u) + \tau(u) + \alpha\langle u,u\rangle_0\right\rbrace  .
\end{equation*}
Let now $E_{k,k=0,1,2...}$ be the collection of eigenspaces of the operator $\Delta_0^V$ with corresponding eigenvalues
$\lambda_{k,k=0,1,2...}$. Again, we denote the eigenspaces and the orthogonal projections onto the eigenspaces by the
same symbol. Then, taking $\eps$ so small that $\eps^2 \leq 1 - \lambda_{0}/\lambda_{1}$, we obtain by
$$
\frac{\lambda_{k}-\lambda_{0}}{\eps^2}=\lambda_k\frac{1-\lambda_{0}/\lambda_{k}}{\eps^2}\geq \lambda_k\frac{1-\lambda_{0}/\lambda_{1}}{\eps^2}\geq \lambda_k
$$
that
\begin{eqnarray*}
\frac{1}{2\eps^2}q_V(u)&=& \sum_{k\geq
1}\frac{\lambda_{k}-\lambda_{0}}{2\eps^2}\langle u, E_{k} u\rangle_{0} \\
&\geq& \frac{1}{2} \sum_{k\geq 1}\lambda_{k}\langle u, E_{k} u\rangle_{0} .
\end{eqnarray*}
That implies by the assumption on $\alpha$ that
\begin{eqnarray*}
\frac{1}{2}\left\lbrace\eps^{-2}q_V(u) + \alpha
\,\langle u,u\rangle_{0} \right\rbrace &\geq& \frac{\lambda_{0}\langle u,E_{0}u\rangle_{0} +\sum_{k\geq
1}\lambda_{k}\langle u, E_{k} u\rangle_{0}}{2} \\
&=& \frac{1}{2}\left(q_0(u) - \tau (u)\right).
\end{eqnarray*}
Thus
\begin{equation*}
F^0_{\eps,\alpha}(u) \geq \frac{1}{2}\left(q_{0}(u) - \tau (u)\right) + \frac{1}{2}\tau (u)= \frac{1}{2}q_{0} (u).
\end{equation*}
\end{Proof}

The next statement is the estimate for the limes inferior of the functionals $F^0_{\eps,\alpha }$ as $\eps$ tends to
zero.

\begin{Lemma}\label{liminf} Let $u_{n,n\geq 1}$ be a weakly convergent sequence in $\bsob^1 (L(1),g_{0})$ with limit
$u$. Let $\eps_{n}\to 0$ be sequence of numbers $0< \eps_{n}\leq 1$. Then $$
\liminf_{n}  F^0_{\eps_{n},\alpha }(u_{n}) \geq F^0_{\alpha} (u), $$
where $$
F^0_{\alpha} (u) = \left\lbrace
\begin{array}{ll}  \frac{1}{2}\int_L \m_L(dp)\left(\langle \d v,\d v\rangle_{L} + \alpha\,v^2\right) & ,u=u_0\overline{v}\in
E_{0}\\
\infty & ,{\mathrm{else}}
\end{array}\right. . $$
\end{Lemma}

\begin{Proof} For $\eps \leq 1 - \lambda_{0}/\lambda_{1}$, we have
$$
\frac{\lambda_{k}-\lambda_{0}}{\eps^2}=\frac{\lambda_k}{\eps}\frac{1-\lambda_{0}/\lambda_{k}}{\eps}\geq \frac{\lambda_k}{\eps}\frac{1-\lambda_{0}/\lambda_{1}}{\eps}\geq \frac{\lambda_k}{\eps}
$$
and therefore
\begin{eqnarray}\label{estimat}
\nonumber\frac{1}{2\eps^2}q_V(u)&=& \sum_{k\geq
1}\frac{\lambda_{k}-\lambda_{0}}{2\eps^2}\langle u, E_{k} u\rangle_{0} \geq \frac{1}{2\eps}\sum_{k\geq
1}\lambda_{k}\langle u, E_{k} u\rangle_{0} \\
&=& \frac{1}{2\eps} \lbrace q_V( u^{\perp}) + \lambda_0 \langle u^{\perp}, u^{\perp}\rangle_0 \rbrace
\end{eqnarray}
where $u^{\perp} = (1-E_0)u$. Now $r_V(u) := q_V( u) + \lambda_0 \langle u, u\rangle_0$, $\tau (u)$ and $\langle
u,u\rangle_{0}$ are all non-negative quadratic forms which are continuous in $\bsob^1 (L(1),g_{0})$. For instance for
the first quadratic form,  $r_V((u_{n}-u)^{\perp})\geq 0$ implies by polarization $$
r_V(u_{n}^{\perp}) \geq  2r_V( u_{n}^{\perp}, u^{\perp}) - r_V(u^{\perp}) = 2 r_V( u_{n}, u^{\perp}) - r_V( u^{\perp})$$
where we denote the associated bilinear form $r_V(-,-)$ by the same symbol. By weak convergence, we have $\lim r_V ( u_{n},  (1-E_0) u) = r_V(u,u^{\perp})$. Now by (\ref{formvertikal})
$$
r_V(u,u^{\perp})=\int_{L}\m_L(dq) \langle d_q u_q, d_q u_q^{\perp}\rangle_q = \int_{L}\m_L(dq) \langle  u_q, D_q u_q^{\perp}\rangle_q .
$$
By decomposability of the projection $E_0$, we have
$$
u^{\perp}_q = (u^{\perp})_q = ((1-E_0)u)_q = u^{\perp}_q-E_{0,q}u_q
$$
which implies $r_V(u,u^{\perp})=r_V(u^{\perp},u^{\perp})=r_V(u^{\perp})$. Thus, we obtain $$
\liminf_{n} r_V (u_{n}^{\perp}) \geq r_V( u^{\perp} ) . $$
Essentially the same argument for the other two quadratic forms yields $$
\liminf_{n}  F^0_{\eps_{n},\alpha}(u_{n}) \geq \frac{1}{2}\left(\tau (u) + \alpha\langle
u,u\rangle_{0} + \liminf_{n} \frac{1}{\eps_{n}} r_V( u_{n}^{\perp}) \right) $$
For $u\in E_0$, by $r_V( u_{n}^{\perp}) \geq 0$ there is nothing more to prove. For
$u\notin E_0$, we have $$
\liminf_{n} r_V ( u_n^{\perp}) = a > 0 $$
and hence for some suitable $a > \delta > 0$ $$
\liminf_{n} \frac{1}{\eps_{n}}r_V( u_{n}^{\perp})   > \liminf_{n} \frac{1}{\eps_{n}}(a
-\delta)  = \infty. $$
It finally remains to show that for $u=u_0\,\overline{v}\in E_0$, $\tau (u) = \langle \d v,\d v\rangle_{L}$. This follows immediately from the local representation of $\tau$ in Lemma \ref{lokalfnull} above and the fact that eigenfunctions $u\in E_0$ are invariant with respect to rotations of the fiber and therefore annihilated by the Killing vector field, i.e. $L_{\mu}^{\nu} u =0$. That implies finally
\begin{eqnarray*}
\tau (u) + \alpha \langle u,u\rangle_0 &=& \int_{L(1)} \m_0(dq) (u_0^2 \overline{\langle \d v,\d v\rangle_L} + \alpha u_0^2 \,\overline{v}^2) \\
&=& \int_{L} \m_L(dp) (\langle \d v,\d v\rangle_L + \alpha  \,v^2)\,\Vert u_0\Vert^2_p \\
&=& \int_{L} \m_L(dp) (\langle \d v,\d v\rangle_L + \alpha  \,v^2).
\end{eqnarray*}
\end{Proof}

\subsection{A Kato type inequality for the induced metric}

By the calculation of the induced metric in local
coordinates, the family $F_{\eps,\alpha}$ of quadratic forms
associated to the rescaling of the induced metric can be decomposed by
\begin{equation}
\label{induquad} F_{\eps,\alpha}(u) = F_{\eps,\alpha}^0 (u) + \frac{1}{2}\left( \langle \d u, \d u\rangle_H +
\eps\,\langle \d u, \d u\rangle_{R(\eps)} + \langle u, W_{\eps} u\rangle_{0}\right)
\end{equation}
where $H$ is given by $H(\eps)= H + \eps \,R(\eps)$ where $H(\eps) = g(\eps) - g_0(\eps)$ from Corollary \ref{Asymptotik
Metrik}, $W_{\eps}$ is the potential in Proposition \ref{modiDir1} and $R(\eps)$ is a symmetric 2-tensor converging
uniformly to some bounded $R(0)$ on $L(1)$ as $\eps$ tends to zero. The central idea to prove convergence of the
minimizers for the functionals associated to the induced metric $g$ is now to compare these functionals with the
functionals associated to the reference metric $g_0$ asymptotically. To control the difference of both functionals as
$\eps$ tends to zero, we will use the following Kato-type inequality. Recall that we write $\overline{W_L}=W_L\circ\pi$
in accordance with the way to indicate basic functions.

\begin{Proposition}\label{Kato} Let $\alpha \geq\lambda_0$. Then there is a constant $K > 0$ such that $$
\left\vert F_{\eps,\alpha}(u)-F^0_{\eps,\alpha}(u)-\frac{1}{2}\langle u,\overline{W_L},u\rangle_{0} \right\vert\leq
K\,\eps\,F^0_{\eps,\alpha}(u). $$
for all $\eps \leq 1 - \lambda_0/\lambda_1$.
\end{Proposition}

To prove this, we establish the relevant inequalities by the two lemmas below. Note that the first statement
$\langle \d u,\d u\rangle_{H} = 0$ for $u\in E_{0}$ which holds due to the fact that functions in $E_0$ are invariant with respect to rotations of the fibers is the crucial one. If the zero order contribution would not vanish on the asymptotic subspace $E_0$, we would not be able to establish the Kato bound above.

\begin{Lemma}\label{KI1} There is a constant $A > 0$ such that $$
\left\vert\langle \d u, \d u\rangle_H \right\vert\leq A\eps
 \,F_{\eps,\alpha}^0 (u) $$
for $\eps \leq 1 - \lambda_0 / \lambda_1$.
\end{Lemma}

\begin{Proof} In local coordinates, we have
$$
\langle \d u,\d u\rangle_H = \frac{1}{3}\int_{L(1)} \m_0(dq)\, w^{\alpha}w^{\beta}R_{\alpha\mu\beta\nu}\partial_{\mu}u\partial_{\nu}u .
$$
By the same idea as for the form $\tau$ in Lemma \ref{lokalfnull}, we use the symmetry properties
$$
R_{\alpha\mu\beta\nu} = - R_{\mu\alpha\beta\nu} = - R_{\alpha\mu\nu\beta}
$$
of the curvature to establish
$$
w^{\alpha}w^{\beta}R_{\alpha\mu\beta\nu}\,\partial_{\mu}u\,\partial_{\nu}u = \frac{1}{4}R_{\alpha\mu\beta\nu}L_{\mu}^{\alpha}u L_{\nu}^{\beta}u
$$
with the Killing vector fields $L_{\mu}^{\alpha} = w^{\alpha}\partial_{\mu} - w^{\mu} \partial_{\alpha}$. By $L_{\mu}^{\alpha}u = 0$ for all $u\in E_0$, that implies
$$
\langle \d u,\d u\rangle_H = \langle \d u^{\perp},\d u^{\perp}\rangle_H .
$$
One way to write the form $\langle \d u,\d u\rangle_H$ in an invariant way uses the fact that the second fundamental form of the embedding of the fibers $\pi^{-1}(q)$ vanishes at the basepoint $q\in L$. Hence, by the {\em Gauss equations}, we can think of $R_{\alpha\mu\beta\nu}$ as the coefficients of the fiber curvature $R_q = R_{\pi^{-1}(q)}$ at $q$. Letting now $X_u$ be the vector field $X_u (q) := \langle \d_q u_q, -\rangle_q$, we obtain at $\xi\in N_qL$
$$
\langle \d u,\d u\rangle_H = \frac{1}{3}\int_{L(1)} \m_0(dp)\,\langle R_{\pi(p)}(\xi, X_u)\xi,X_u\rangle .
$$
That implies
\begin{eqnarray*}
\vert \langle \d u,\d u\rangle_H \vert &=& \frac{1}{3}\int_{L(1)} \m_0(dp)\,\langle R_{\pi(p)}(\xi, X_u)\xi,X_u\rangle \\
&\leq& A \,\int_{L} \m_L(dq)\,\Vert X_u\Vert_q^2 \\
&=& A \,\int_{L} \m_L(dq)\,\Vert \d_q u_q\Vert_q^2 \\
&=& A\, r_V(u) .
\end{eqnarray*}
Now, by combining this estimate with(\ref{estimat}), we obtain
$$
\langle \d u,\d u\rangle_H = \langle \d u^{\perp},\d u^{\perp}\rangle_H\leq A \,r_V(u^{\perp}) \leq A\eps\frac{1}{\eps^2} q_V(u) \leq A\,\eps\, F_{\eps,\alpha}^{0}(u).
$$
\end{Proof}

The second lemma is rather straightforward using no special knowledge on the structure of the functionals under
consideration.

\begin{Lemma}\label{KI2} There are constants $B,C > 0$ such that for $\eps^2 \leq 1-\lambda_0/\lambda_1$ and $\alpha
\geq \lambda_0$ we have
\begin{enumerate}
\item $\left\vert\langle \d u, \d u\rangle_{R(\eps)}\right\vert\leq B
\,F^0_{\eps,\alpha}(u)$, \item $\left\vert\langle u, (W_{\eps}-\overline{W_L}) u\rangle_{0}\right\vert\leq C\eps\,
F^0_{\eps,0}(u)$.
\end{enumerate}
\end{Lemma}

\begin{Proof} {\rm (i)} The first inequality follows from the fact that $R(\eps)$ converges uniformly to $R(0)$ and thus by Lemma \ref{koerziv}
$$
\left\vert\langle \d u, \d u\rangle_{R(\eps)}\right\vert\leq  2B\,
q_0(u)\leq B \,F^0_{\eps,\alpha}(u). $$
{\rm (ii)} By uniform convergence of $W_{\eps}$ to $W_0$ as $\eps$ tends to zero, Taylor expansion in the normal Fermi coordinates yields by the fact that the Sobolev norm on $\bsob^1(L(1),g_0)$ dominates the $L^2$-norm and can be equivalently described by $q_0$ ({\em Poincare - inequality})
$$
\left\vert\langle u, (W_{\eps}-\overline{W_L}) u\rangle_{0}\right\vert\leq D''\eps\,\Vert u\Vert_{0}^2 \leq D'
\eps \, q_0 (u) \leq D \eps F^0_{\eps,0}(u). $$
\end{Proof}

The proof of Proposition \ref{Kato} is now

\begin{Proof} By Lemma \ref{KI1} and \ref{KI2}, we have for $\eps >0$ small enough
\begin{eqnarray*}
& & \left\vert F_{\eps,\alpha}(u) - F_{\eps,\alpha}^0 (u) - \frac{1}{2}\langle u, \overline{W_L} u\rangle_{0}
\right\vert \\
&\leq& \frac{1}{2}\left\lbrace\left\vert\langle \d u, \d u\rangle_H\right\vert + \eps\,\left\vert\langle
\d u, \d u\rangle_{R(\eps)}\right\vert + \left\vert\langle u, (W_{\eps}-\overline{W_L})
u\rangle_{0}\right\vert\right\rbrace \\
&\leq& \frac{1}{2}\left\lbrace A \eps \,F_{\eps,\alpha}^0 (u) + B
\eps\,F^0_{\eps,\alpha}(u)+C\eps\, F^0_{\eps,0}(u)\right\rbrace \\
&\leq& K\,\eps\,F^0_{\eps,0}(u).
\end{eqnarray*}
\end{Proof}

\subsection{Equi-coercivity and the proof of Theorem \ref{MainMosco}} In this subsection, we prove Theorem \ref{MainMosco} in the form of Proposition \ref{epiMAIN} below and equi-coercivity
of the corresponding sequence. These are the facts that are necessary to conclude convergence of the associated
semigroups in the next section. We begin by proving Proposition \ref{koerzmain}, equi-coercivity of the sequence of functionals $F_{\eps,\alpha}$ in the following slightly stronger form.

\begin{Lemma}\label{kroz2} Let $\alpha \geq \lambda_0 +
\sup_{q\in L} \lbrace - W(q)\wedge 0\rbrace$. Then $$
F_{\eps,\alpha}(u)\geq \frac{1}{4}q_{0}(u) $$
for all $\eps \leq 1/ 2K$.
\end{Lemma}

\begin{Proof} Let $\alpha = \alpha_1 + \alpha_2$ with $\alpha_1 \geq \lambda_0$ and $\alpha_2 \geq \sup_{q\in L} \lbrace
- W(q)\wedge 0\rbrace$. By Proposition \ref{Kato}, we have
\begin{eqnarray*}
& & F_{\eps,\alpha}(u) \\
&=& F^0_{\eps,\alpha_1}(u) + F_{\eps,\alpha_1}(u)-F^0_{\eps,\alpha_1}(u)-\frac{1}{2}\langle u,\overline{W_L}
u\rangle_{g_0} +\frac{1}{2}\langle u,(\overline{W_L}+\alpha_2)u\rangle_{g_0}\\
&\geq& (1 - K\,\eps)\,F^0_{\eps,\alpha_1}(u)+\frac{1}{2}\langle u,(\overline{W_L}+\alpha_2) u\rangle_{g_0}\\
&\geq& (1 - K\,\eps)\,F^0_{\eps,\alpha_1}(u)
\end{eqnarray*}
and $1 - K\eps \geq 1/2$ for $\eps \leq 1/ 2K$. By Lemma
\ref{koerziv}, that implies the statement.
\end{Proof}

Now we prove condition (ii) from Definition \ref{epidef} with an explicit bound for $\alpha_0$.

\begin{Lemma}\label{liminfMAIN} Let $\alpha \geq \lambda_0 +
\sup_{q\in L} \lbrace - W(q)\wedge 0\rbrace$ and $u_{n,n\geq 1}$ be a weakly convergent sequence in $\bsob^1
(L(1),g_{0})$. Let $\eps_{n}\to 0$ be a sequence of numbers $0< \eps_{n}\leq 1$. Then $$
\liminf_{n} F_{\eps_n,\alpha}(u_{n}) \geq F_{\alpha} (u), $$
where $$
F_{\alpha} (u) = \left\lbrace
\begin{array}{ll}
\frac{1}{2}\int_L \m_L(dq)\left(\langle \d v,\d v\rangle_{L} + (\overline{W_L}  + \alpha)v^2\right) & ,u=u_0\overline{v}\in E_{0}\\
\infty & ,{\mathrm{else}}
\end{array}\right. . $$
\end{Lemma}

\begin{Proof} Let $\alpha = \alpha_1 + \alpha_2$ with $\alpha_1 \geq \lambda_0$ and $\alpha_2 \geq \sup_{q\in L} \lbrace
-W(q)\wedge 0\rbrace$. By Proposition \ref{Kato} and Lemma \ref{liminf} we have
\begin{eqnarray*}
& & \liminf_{n} F_{\eps_n,\alpha }(u) \\
&=& \liminf_{n} F_{\eps_n,\alpha_1}(u)-\frac{1}{2}\langle u_n,\overline{W_L} u_n\rangle_{0} + \frac{1}{2}\langle
u_n,(\overline{W_L}+\alpha_2) u_n\rangle_{0}\\
&\geq & \liminf_{n}(1 - K\,\eps_n)\,F^0_{\eps_n,\alpha_1}(u) +
\liminf_n \frac{1}{2}\langle u_n,(\overline{W_L}+\alpha_2) u_n\rangle_{0}\\
&\geq& F^0_{\alpha_1} (u) + \frac{1}{2}\langle u,(\overline{W_L}+\alpha_2) u\rangle_{0} \\
&=& F_{\alpha} (u),
\end{eqnarray*}
using
$$
\langle u,\overline{W_L}u\rangle_{0}=\int_L \m(dp) W_L\, v^2\,\Vert u_0\Vert_p^2 = \int_L \m(dp) W_L\, v^2.
$$
for the last step.
\end{Proof}

The following proposition is the reformulation of Theorem \ref{MainMosco} with the bound for $\alpha_{0}$ from Lemma \ref{liminfMAIN}.

\begin{Proposition}\label{epiMAIN} Let $\alpha \geq \alpha_{0}:=\lambda_0 +
\sup_{q\in L} \lbrace - W(q)\wedge 0\rbrace$. Then the functionals $F_{\eps,\alpha}$ epi-converge to $F_{\alpha}$ with
respect to the weak topology on $\bsob^1(L(1),g_0)$.
\end{Proposition}

\begin{Proof} First of all, by Proposition \ref{koerzmain} or likewise Lemma \ref{kroz2}, the functionals are uniformly bounded below by the 1-Sobolev norm. Since $\bsob^1 (L(1),g_0)$ is a reflexive Banach space, norm-bounded sets are relatively compact in the weak topology by
the Banach - Alaoglu theorem. Hence, by \cite{DaM:92}, Proposition 7.7, p.70, the sequence $F_{\eps_{n},\alpha}$ is
equi-coercive. By the preceding definition it remains to prove:

\vspace{0.1cm}

{\rm (i)} For all $u\in\bsob^1(L(1),g_0)$ there is a weakly convergent sequence $u_n\to u$ such that
\begin{equation}
\lim_{n}  F_{\eps_{n},\alpha}(u_{n}) = F_{\alpha} (u).
\end{equation}
This follows by pointwise convergence of the functionals. Hence we may take $u_n\equiv u$.

{\rm (ii)} For all $u\in\bsob^1(L(1),g_0)$ and for all weakly convergent sequences $u_n\to u$ we have
\begin{equation*}
\liminf_{n}  F_{\eps_{n},\alpha }(u_{n}) \geq F_{\alpha} (u).
\end{equation*}
This is the statement of Lemma \ref{liminfMAIN}.
\end{Proof}

\subsection{Proof of Theorem \ref{MainSemi}}

\begin{Proof} {\rm (i)} First of all, the right hand side of the equation is meant to be $$
E_0 \,(\Delta_L + W_L + \alpha)^{-1}\, E_0 u = u_0 \,\left\lbrace(\Delta_L + W_L + \alpha)^{-1}f\right\rbrace\circ\pi $$
where we use that $E_0u = u_0\,\overline{f}$ with $f\in L^2 (L(1),g_0)$. By
Proposition \ref{koerzmain} and the remark following Definition \ref{epidef}, the functionals $$
F_{\eps,\alpha }(u)- \langle w,u\rangle_{0} \geq \frac{1}{8} q_0( u ) - 2\Vert  w
\Vert_{\sob^1(L(1),g_0)^*}^2 $$
where $\sob^1(L(1),g_0)^*$ denotes the dual Sobolev space, are uniformly bounded below. That implies that for all
$r\in\R$ the sets $$
K_r:=\left\lbrace F_{\eps,\alpha }(u)- \langle w,u\rangle_{0} \leq r\right\rbrace \subset\subset \bsob^1 (L(1),g_0)
$$
are norm bounded and therefore relatively compact in the weak topology of the reflexive Banach space
$\bsob^1(L(1),g_0)$. That means, the sequence is {\em equi-coercive}. By Proposition
\ref{epiMAIN}, again together with the remark above, the functionals $F_{\eps,\alpha,w}:=F_{\eps,\alpha }-\langle
w,-\rangle_{0}$ epi-converge to $F_{\alpha,w} := F_{\alpha}- \langle w,-\rangle_{0}$ with respect to the weak topology
on $\bsob^1(L(1),g_0)$.

To compute the minimizers of $F_{\eps,\alpha,w}$, note that these functions are strictly convex and
differentiable with a gradient
$$
\nabla F_{\eps,\alpha,w} : \bsob^1 (L(1),g_0)\rightarrow \bsob^1 (L(1),g_0)^* $$
which is at $u\in\bsob^1(L(1),g_0)$ given by
$$
\nabla_u F_{\eps,\alpha,w}\lbrack x \rbrack = \eps^{-2}\, q_V (u,x)+
\tau (u,x)  + \langle u, (W_L+\alpha)x \rangle_{0} -\langle w,x \rangle_{0}.
$$
Hence the condition for $u^*_{\eps}$ to be a minimizer (which by strict convexity implies that it is in fact {\em the}
minimizer of the functional), is given by $\nabla_{u^*_{\eps}}F_{\eps,\alpha,w}
\equiv 0$.

\noindent This is equivalent to $\nabla_{u^*_{\eps}}F_{\eps,\alpha,w}\lbrack x \rbrack = 0$ for all $x\in\bsob^1
(L(1),g_0)$ which implies that $u^*_{\eps}$ is a {\em weak solution} of $$
(\Delta (\eps) + \alpha ) u^*_{\eps} = w. $$
Since $\Delta(\eps)$ is basically the Laplace-Beltrami operator associated to the metric $g(\eps)$, $\Delta (\eps) +
\alpha $ is {\em elliptic} and the solution is indeed a strong solution. Furthermore, by Friedrich's construction,
$\Delta (\eps)+ \alpha$ is self-adjoint and positive on $L^2 (L(1),g_0)$. Hence, it is invertible and the minimizer is
given by $u^*_{\eps}=(\Delta (\eps) + \alpha )^{-1} w$.

\noindent The limit functional $F_{\alpha,w}$ is strictly proper convex, i.e. strictly convex on the subspace $E_0$
where it is finite. Hence there will be one minimizer $u^*_0\in E_0$. By Corollary \ref{kern Dirichlet}, a function
$u\in E_0$ can be written $u=u_0\,\overline{f}$ and hence we may write equivalently using $E_0w = u_0\,\overline{g}$
$$
F_{\alpha,w}(u) = F_{\alpha,E_0w}(u)=\frac{1}{2}\int_L\m_L(dp)\left(\langle \d f,\d f\rangle_{L} + (\alpha + W_L)f^2\right) - \langle g,f\rangle_{L} .
$$
As gradient on the subspace $E_0$ we therefore obtain
$$
\nabla_u F^0_{\alpha,w}\lbrack v \rbrack = \int_L\m_L(dp)\left(\langle \d f,\d h\rangle_{L} + (\alpha + W_L)fh\right) - \langle g,h\rangle_{L}
$$
where $v= u_0\, \overline{h}\in E_0$. Hence, the minimizer of the limit functional is given by $u_0^* = u_0\, \overline{f^*}$
where by the same arguments as above, $f^*$ is a weak solution of $$
(\Delta_L + W_L + \alpha) f^* = g. $$
By ellipticity of the Laplace - Beltrami operator on a closed manifold, this is again a strong solution and since the
Laplace - Beltrami operator of a closed manifold is self - adjoint and non-negative, $\Delta_L + \alpha$ is invertible
on $L^2 (L(1), g_0)$ and we obtain $u_0^* = u_0\,(\Delta_L + \alpha)^{-1}g$ or equivalently
$$
u_0^* =E_0(\Delta_L + \alpha)^{-1}\,E_0 w .
$$
In particular, the minimizer of the limit functional is unique and by \cite{DaM:92}, Corollary 7.24, p. 84, each
sequence of minimizers $u_{\eps_n}^*$ of $F^0_{\eps_n,\alpha,w}$ converges to $u_0^*$ as $n$ tends to infinity with
respect to the weak topology of $\bsob^1(L(1),g_0)$. By \cite{DaM:92}, Corollary 8.8, p. 92 (cf. in particular Example
8.9 on the same page), this is equivalent to norm-convergence on $L^2(L(1),g_0)$ since on bounded sets, the weak
topology on $\bsob^1(L(1),g_0)$ coincides with the one induced by the $L^2$-norm. That implies the first statement.

\noindent{\rm (ii)} By Proposition \ref{koerzmain}, the associated family of quadratic forms is non-negative and thus the family of self-adjoint operators $\Delta (\eps)+\alpha_{0}$ is uniformly bounded below, too. Thus the spectra $\spec{\Delta(\eps)+\alpha_{0}}$
are all contained in $\lbrack 0,\infty)$. That implies that the operators $-(\Delta(\eps)+\alpha_{0})$ are {\em
sectorial} in the sense of \cite{EngNag:00}, 4.1 Definition, p. 96 and that they are even {\em uniformly sectorial}
meaning that:
\begin{enumerate}
\item There is a common sector ($0<\delta \leq\pi/2$) $$
\Sigma_{\pi/2 +\delta}:=\lbrace z\in\C\,:\,\vert \arg (z)\vert < \pi/2 + \delta\rbrace - \lbrace 0 \rbrace $$
which is contained in all the resolvent sets $\rho (-(\Delta(\eps)+\alpha_{0}))$, $\eps > 0$. (In our case we can take
for instance $\delta = \pi /4$.)
\item For all $\eta\in (0,\delta)$ there is some $M_{{\eta}}>0$ such that $$
\Vert R(-(\Delta(\eps)+\alpha_{0}),-z)\Vert \leq \frac{M_{\eta}}{\vert z \vert} $$
for all $0\neq z\in\overline{\Sigma_{\pi/2 +\delta-\eta}}$.
\end{enumerate} To prove the second assertion, note that by $$
\Vert R(-(\Delta(\eps)+\alpha_{0}),-z)\Vert \leq \frac{\sqrt{2}}{\vert z \vert} $$
independent of $\eps > 0$, all these operators generate {\em analytic
semigroups} in the sense of \cite{EngNag:00}, Ch. II, 4a, p. 96 ff. by the {\em Dunford-integral} $$
e^{-\frac{t}{2}(\Delta(\eps)+\alpha_{0})}u:=\frac{1}{2\pi i}\int_{\gamma}dz
\,e^{\frac{tz}{2}}R(-(\Delta(\eps)+\alpha_{0}),-z)u $$
and $$
E_{0}\,e^{-\frac{t}{2}(\Delta_L + W_L+\alpha_{0})}\,E_{0}u:=\frac{1}{2\pi i}\int_{\gamma}dz
\,E_{0}\,e^{\frac{tz}{2}}R(-(\Delta_L + W_L+\alpha_{0}),-z)\,E_{0}u, $$
respectively. Here, $\gamma$ denotes a suitable curve contained in the uniform sector $\Sigma_{3\pi/4}$ independent of $\eps > 0$ (cf.
\cite{EngNag:00}, 4.2 Definition, p.96). Since the resolvents are holomorphic functions on the common sector, the convergence result (i) for
the resolvents implies that for all $u\in L^2(L(1),g_{0})$ $$
\lim_{\eps\to 0}R(-(\Delta(\eps)+\alpha_{0}),-z)u = E_{0}\,R(-(\Delta_L + W_L+\alpha_{0}),-z)\,E_{0}u $$
for all $z\in \Sigma_{3\pi/4}$. By the choice of the curve and by the uniform estimate of the resolvent norm above, we hence obtain using dominated convergence
$$
\lim_{\eps\to 0}\frac{1}{2\pi i}\int_{\gamma}dz\,e^{\frac{tz}{2}}(R(-(\Delta(\eps)+\alpha_{0}),-z) - E_{0}\,R(-(\Delta_L
+ W_L+\alpha_{0}),-z)\,E_{0})u = 0
$$
which yields (ii) after division by $e^{-\alpha_{0}t/s}$.

\noindent{\rm (iii)} As said above, the operators $\Delta (\eps)$ are uniformly bounded below. Hence, the operator norms of the semigroups
$$\Vert e^{-\frac{t}{2}\Delta (\eps)}\Vert < C $$
are uniformly bounded above and thus
$$
\Vert e^{-\frac{t}{2}\Delta (\eps)}u_{\eps} - e^{-\frac{t}{2}\Delta (\eps)}u\Vert < c\,\Vert u_{\eps} - u\Vert \to 0
$$
as $\eps\to 0$.
\end{Proof}

\section{Acknowledgement.} The author wishes to thank Martin Kolb and Jakob Wachsmuth for pointing out substantial errors in earlier versions of the paper.

\bibliographystyle{plain}

\end{document}